\documentclass[13pt]{article}
\usepackage{latexsym}
\usepackage{geometry}
\usepackage{graphicx}
\usepackage{amsmath, amssymb, amsthm}
\usepackage{booktabs}
\usepackage{algorithm}
\usepackage{algorithmic}

\usepackage[utf8]{inputenc} 
\usepackage[T1]{fontenc}    

\usepackage{url}
\usepackage{natbib}

\usepackage{appendix}

\usepackage{amsfonts}
\usepackage{multirow}
\usepackage{multicol}

\usepackage{color}
\usepackage{algorithm}
\usepackage{algorithmic}

\usepackage{xcolor,colortbl}
\definecolor{LightCyan}{rgb}{0.88,1,1}

\usepackage{graphicx}
\usepackage{subfig}
\usepackage{hyperref}
\usepackage{tcolorbox}

\newtheorem{theorem}{Theorem}
\newtheorem{lemma}{Lemma}

\newtheorem{definition}{Definition}
\newtheorem{assumption}{Assumption}
\newtheorem{remark}{Remark}

\begin{document}
\title{ Adaptive Mirror Descent Bilevel Optimization}
\author{Feihu Huang\thanks{Feihu Huang is with College of Computer Science and Technology,
Nanjing University of Aeronautics and Astronautics, Nanjing, China;
and also with MIIT Key Laboratory of Pattern Analysis and Machine Intelligence, Nanjing, China. Email: huangfeihu2018@gmail.com} }

\date{ }


\maketitle

\begin{abstract}
In the paper, we propose a class of efficient adaptive bilevel methods based on mirror descent for nonconvex bilevel optimization, where its upper-level problem is nonconvex possibly with nonsmooth regularization, and its lower-level problem is also nonconvex while satisfies Polyak-{\L}ojasiewicz (PL) condition.
To solve these deterministic bilevel problems, we present an efficient adaptive projection-aid gradient (i.e., AdaPAG) method based on mirror descent, and prove that it obtains the best known gradient complexity of $O(\epsilon^{-1})$ for finding an $\epsilon$-stationary solution of nonconvex bilevel problems.
To solve these stochastic bilevel problems, we propose an efficient adaptive stochastic projection-aid gradient (i.e., AdaVSPAG) methods based on mirror descent and variance-reduced techniques, and prove that it obtains  the best known gradient complexity of $O(\epsilon^{-3/2})$ for finding an $\epsilon$-stationary solution.
Since the PL condition relaxes the strongly convex, our algorithms can be used to nonconvex strongly-convex bilevel optimization.
Theoretically, we provide a useful convergence analysis framework for our methods under some mild conditions, and prove that our methods have a fast convergence rate of $O(\frac{1}{T})$, where $T$ denotes the number of iterations.
\end{abstract}

\section{Introduction}
Bilevel optimization~\citep{colson2007overview} is an useful two-level hierarchical optimization,
is widely used in many machine learning tasks such as hyperparameter learning~\citep{franceschi2018bilevel},
 meta learning~\citep{franceschi2018bilevel,ji2021bilevel} and reinforcement learning~\citep{hong2020two,chakraborty2023aligning}.
In the paper, we consider a class of nonconvex bilevel optimization problem, defined as
\begin{align}
& \min_{x \in \mathbb{R}^d, y\in y^*(x)} \ f(x,y) + h(x),  & \mbox{(Upper-Level)} \label{eq:1} \\
& \mbox{s.t.} \ y^*(x) \equiv \arg\min_{y\in \mathbb{R}^p} \ g(x,y),  & \mbox{(Lower-Level)} \nonumber
\end{align}
where $f(x,y)$ with $y\in y^*(x)$ is possibly nonconvex, and $h(x)$ is a convex but possibly nonsmooth regularization such as $h(x)=\|x\|_1$ or $h(x)=0$ when $x \in \mathcal{X}\subseteq \mathbb{R}^d$ with convex set $\mathcal{X}$ otherwise $h(x)=+\infty$. Here the lower function $g(x,y)$ is possibly nonconvex on any $y$ but  satisfies Polyak-Lojasiewicz (PL) condition~\citep{polyak1963gradient}, which
 relaxes the strong convexity and is used some machine learning models such the over-parameterized deep neural networks~\citep{frei2021proxy,song2021subquadratic}
have been shown to satisfy the PL condition. Meanwhile, we also study the stochastic version of Problem~(\ref{eq:1}):
\begin{align}
& \min_{x \in \mathbb{R}^d, y\in y^*(x)} \ \mathbb{E}_{\xi}\big[f(x,y;\xi)\big] + h(x),  & \mbox{(Upper-Level)} \label{eq:2} \\
& \mbox{s.t.} \ y^*(x) \equiv \arg\min_{y\in \mathbb{R}^p} \ \mathbb{E}_{\zeta}\big[g(x,y;\zeta)\big],  & \mbox{(Lower-Level)} \nonumber
\end{align}
where $f(x,y)=\mathbb{E}_{\xi}\big[f(x,y;\xi)\big]$ is possibly nonconvex, and $h(x)$ is a convex but possibly nonsmooth regularization, and $g(x,y)=\mathbb{E}_{\zeta}\big[g(x,y;\zeta)\big]$ is nonconvex on any $y$ and satisfies PL inequality. Here $\xi$ and $\zeta$ are independently random variables following unknown distributions $\mathcal{D}$ and $\mathcal{S}$, respectively.
In fact, these Problems~(\ref{eq:1}) and~(\ref{eq:2}) are widely appeared in many machine learning problems
such as few-shot mate learning~\citep{huang2023momentum} and policy learning~\citep{chakraborty2023aligning}.

\begin{table*}
  \centering
  \caption{ Comparison of gradient \textbf{complexity }
between our algorithms and the existing (adaptive) bilevel algorithms in solving nonconvex-PL bilevel problems in finding an $\epsilon$-stationary solution ($\|\mathcal{G}(x,\nabla F(x),\gamma)\|^2\leq \epsilon$ or its equivalent variants, where $F(x)=f(x,y)$ with $y\in y^*(x)$). Here $g(x,\cdot)$ denotes the function on the second variable $y$ with fixing variable $x$. \textbf{Note that} the GALET~\citep{xiao2023generalized} method simultaneously uses the PL condition, its Assumption 2 (i.e., let $\sigma_g = \inf_{x,y}\{\sigma_{\min}^{+}(\nabla^2_{yy} g(x,y))\} >0$ for all $(x,y)$) and its Assumption 1 (i.e., $\nabla^2_{yy} g(x,y)$ is Lipschitz continuous). Clearly, when Hessian matrix $\nabla^2_{yy} g(x,y)$ is singular, its Assumption 1 and Assumption 2 imply that the lower bound of the non-zero singular values $\sigma_g$ is close to zero (i.e., $\sigma_g\rightarrow 0$), under this case, the convergence results of the GALET are \textbf{meaningless}, e.g., the constant $L_w = \frac{\ell_{f,1}}{\sigma_g}+\frac{\sqrt{2}\ell_{g,2}\ell_{f,0}}{\sigma_g^2}\rightarrow + \infty$ used in its Lemmas 6 and 9. Under the other case, the PL condition, Lipschitz continuous of Hessian and its Assumption 2 (the singular values of Hessian is bounded away from 0, i.e., $\sigma_g>0$) imply that the GALET assumes the strongly convex assumption (More detail, please see the following Section 2.1).  \textbf{SC} stands for strongly convex. \textbf{ALR} denotes adaptive learning rate. }
  \label{tab:1}
   \resizebox{\textwidth}{!}{
\begin{tabular}{c|c|c|c|c|c}
  \hline
  \textbf{Problem} & \textbf{Algorithm} & \textbf{Reference} & \textbf{Assumption} $g(x,\cdot)$  & \textbf{Complexity} & \textbf{ALR} \\ \hline \hline
  \multirow{5}*{\textbf{Deterministic}} &  BOME  & \cite{liu2022bome} & PL / local-PL  & $O(\epsilon^{-3/2})$ / $O(\epsilon^{-2})$ &  \\   \cline{2-6}
  & V-PBGD  & \cite{shen2023penalty} & PL / local-PL & $O(\epsilon^{-3/2})$ / $O(\epsilon^{-3/2})$ & \\  \cline{2-6}
  & GALET  & \cite{xiao2023generalized} & SC / PL  & $O(\epsilon^{-1})$ / \textbf{Meaningless} &  \\ \cline{2-6}
  & MGBiO  & \cite{huang2023momentum} & PL / local-PL  &  {\color{red}{$O(\epsilon^{-1})$}} / {\color{red}{$O(\epsilon^{-1})$}} &  \\ \cline{2-6}
& AdaPAG  & Ours  & PL / local-PL & {\color{red}{$O(\epsilon^{-1})$}} / {\color{red}{$O(\epsilon^{-1})$}} & $\surd$  \\  \hline \hline
\multirow{5}*{\textbf{Stochastic}} & V-PBSGD & \cite{shen2023penalty}  & PL & $\tilde{O}(\epsilon^{-3})$ & \\  \cline{2-6}
 & MSGBiO  & \cite{huang2023momentum}  & PL & $\tilde{O}(\epsilon^{-2})$ & \\ \cline{2-6}
 & VR-MSGBiO  & \cite{huang2023momentum}  & PL  & {\color{red}{$\tilde{O}(\epsilon^{-3/2})$}} & \\  \cline{2-6}
 & VR-BiAdam  & \cite{huang2021biadam}  & SC  & {\color{red}{$\tilde{O}(\epsilon^{-3/2})$}} & $\surd$  \\  \cline{2-6}
 & AdaVSPAG  & Ours  & PL  & {\color{red}{$O(\epsilon^{-3/2})$}} & $\surd$ \\  \hline
\end{tabular}
 }
\end{table*}

Given $g\big(x,y\big)=0$, the problem~(\ref{eq:1}) will reduce to the standard single-level
 optimization problem.
Compared with the standard single-level optimization, the major difficulty of bilevel optimization
is that the minimization of the upper and lower-level objectives are intertwined via the
mapping $y^*(x)\in \arg\min_y g(x,y)$.
To handle this difficulty, recently many bilevel optimization methods~\citep{ghadimi2018approximation,hong2020two,ji2021bilevel,chen2022single}
 have been proposed by imposing the
strong convexity assumption on the lower-level objective.
The basic idea of these methods builds on the hyper-gradient,
\begin{align} \label{eq:g1}
 \nabla F(x_t)  = \nabla_xf(x_t,y^*(x_t)) - \nabla^2_{xy}g(x_t,y^*(x_t)) \big[\nabla^2_{yy}g(x_t,y^*(x_t))\big]^{-1}\nabla_yf(x_t,y^*(x_t)).
\end{align}
Thus, there exist different strategies to estimate this hyper-gradient $\nabla F(x_t)$.
Recently, they basically can be divided into two categories:

The first category is to directly estimate the second-order
derivatives of $g(x,y)$. For example, some existing methods~\citep{chen2022single,dagreou2022framework} require an explicit extraction of second-order information
of $g(x,y)$ with a major focus on efficiently estimating its Jacobian and inverse Hessian.
Meanwhile, some methods~\citep{chen2022single,li2022fully,dagreou2022framework} avoid directly estimating its second-order computation and only use the first-order information of both upper and lower objectives.
The second category is to
estimate the optimal solution $y^*(x)$ of the lower-level problem for every given $x$, which updates the variable $y$ multiple times before updating $x$. For example, several methods~\citep{ghadimi2018approximation,hong2020two,chen2021closing,ji2021bilevel} have been developed to
effectively track $y^*(x)$ without waiting for many inner iterations before updating variable $x$.
To solve stochastic bilevel problem~\eqref{eq:2}, some efficient stochastic gradient methods~\citep{ghadimi2018approximation,hong2020two,ji2021bilevel,chen2022single} have been proposed.
Subsequently, \cite{yang2021provably,khanduri2021near} proposed some accelerated stochastic bilevel methods
based on variance reduced techniques. \cite{huang2022enhanced,chen2023accelerated} studied the nonconvex strongly-convex bilevel optimization problems with nonsmooth regularization. Meanwhile, \cite{huang2021biadam} designed a class of efficient adaptive bilevel methods for the nonconvex strongly-convex bilevel optimization.

The above bilevel optimization methods mainly rely on the
restrictive strong convexity condition of lower-level problem. Recently, some bilevel approaches~\citep{liu2020generic,liu2021towards,sow2022constrained,arbel2022non,liu2022bome,
chen2023bilevel,liu2023averaged,lu2023first,shen2023penalty} have been proposed for bilevel optimization without lower-level strong convexity condition. For example, \cite{sow2022constrained} studied nonconvex-convex bilevel optimization with with multiple inner minima. More recently, \cite{liu2023averaged} developed an effective averaged method of multipliers for bilevel optimization with convex lower-level. Meanwhile, \cite{lu2023first} studied the nonconvex constrained bilevel optimization with convex lower-level with nonsmooth regularization.
Moreover, some methods~\citep{liu2022bome,
chen2023bilevel,liu2023averaged,huang2023momentum,kwon2023penalty} further studied the bilevel optimization with non-convex lower-level.
For example,
\cite{liu2022bome} proposed an effective first-order method for nonconvex-PL bilevel optimization, where the lower-level problem is nonconex but satisfies PL condition. Meanwhile, \cite{shen2023penalty} designed an effective penalty-based gradient method for the constrained nonconvex-PL bilevel optimization.
More recently, \cite{huang2023momentum} proposed a class of efficient momentum-based gradient methods for the nonconvex-PL bilevel optimization.

Although the above bilevel methods studied bilevel optimization problems without strong-convexity lower-level, the adaptive bilevel optimization is missing under the non-convexity condition. Meanwhile, the existing adaptive methods~\citep{huang2021biadam} for nonconvex strongly-convex bilevel optimization rely on the global adaptive learning rates, i.e., adaptive matrices $B_t=b_tI_p \ (b_t>0)$ for all $t\geq1$. Thus, there exists a natural question:
\begin{center}
\begin{tcolorbox}
\textbf{ Could we design efficient adaptive bilevel optimization methods without relying on strong-convexity condition and some specific adaptive learning rates? }
\end{tcolorbox}
\end{center}

In the paper, we affirmatively answer to the above question, and propose
a class of efficient adaptive mirror descent bilevel optimization methods
to solve the nonconvex bilevel problems~(\ref{eq:1}) and (\ref{eq:2}) based on the mirror descent~\citep{beck2003mirror}.
Our main contributions are given:
\begin{itemize}
\item[1)] We propose an efficient adaptive projection-aid gradient (i.e.,AdaPAG) method  based on mirror descent to solve the deterministic Problem~(\ref{eq:1}). Since PL condition
 relaxes the Strong Convexity (SC), our AdaPAG method clearly can be used to solve the nonconvex-SC bilevel problems.
 \item[2)] We propose an efficient adaptive stochastic projection-aid gradient (i.e., AdaVSPAG) method based on mirror descent and variance-reduced techniques to solve the stochastic Problem~(\ref{eq:2}).
\item[3)] We provide a solid convergence analysis framework for our methods. Under some mild conditions, we prove that our AdaPAG method reaches the best known gradient complexity of $O(\epsilon^{-1})$ for finding an $\epsilon$-stationary solution of Problem~(\ref{eq:1}). Meanwhile, our AdaVSPAG method obtains the near-optimal gradient complexity of $O(\epsilon^{-3/2})$ for finding an $\epsilon$-stationary solution of Problem~(\ref{eq:2}).
\end{itemize}

\subsection*{Notations}
Matrix $A\succ 0$ is a positive definite. $\|\cdot\|$ denotes the $\ell_2$ norm for vectors and spectral norm for matrices. $\langle x,y\rangle$ denotes the inner product of two vectors $x$ and $y$. For vectors $x$ and $y$, $x^r \ (r>0)$ denotes the element-wise power operation, $x/y$ denotes the element-wise division and $\max(x,y)$ denotes the element-wise maximum. $I_{d}$ denotes a $d$-dimensional identity matrix.
Given function $f(x,y)$, $f(x,\cdot)$ denotes  function \emph{w.r.t.} the second variable with fixing $x$,
and $f(\cdot,y)$ denotes function \emph{w.r.t.} the first variable
with fixing $y$. $a_t=O(b_t)$ denotes that $a_t \leq c b_t$ for some constant $c>0$. The notation $\tilde{O}(\cdot)$ hides logarithmic terms. $\nabla_x$ denotes the partial derivative on variable $x$. For notational simplicity, let $\nabla^2_{xy}=\nabla_x\nabla_y$ and $\nabla^2_{yy}=\nabla_y\nabla_y$. $\Pi_C[\cdot]$ denotes Euclidean projection onto the set $\{x\in \mathbb{R}^d: \|x\| \leq C\}$ for some constant $C>0$.
$\hat{\Pi}_{C}[\cdot]$ denotes a
projection on the set $\{X\in \mathbb{R}^{d\times p}: \|X\| \leq C\}$ for some constant $C>0$.
$\mathcal{S}_{[\mu,L_g]}$ denotes a
projection on the set $\{X\in \mathbb{R}^{d\times d}: \mu \leq \varrho(X) \leq L_g\}$,
where $\varrho(\cdot)$ denotes the eigenvalue function. Both $\hat{\Pi}_{C}$ and $\mathcal{S}_{[\mu,L_g]}$ can be implemented by using Singular Value Decomposition (SVD) and thresholding the singular values.

\section{ Preliminaries }
In this section, we provide some mild assumptions and useful lemmas on the above bilevel problems~(\ref{eq:1}) and~(\ref{eq:2}).

\subsection{ Mild Assumptions}

\begin{assumption} \label{ass:1}
 The function $g(x,\cdot)$ satisfies the Polyak-Lojasiewicz (PL) condition,
 if there exists $\mu > 0$ such that for any given $x$, it holds that
\begin{align}
 \|\nabla_y g(x,y)\|^2 \geq 2\mu \big(g(x,y)-\min_y g(x,y)\big), \ \forall y\in \mathbb{R}^p.
\end{align}
\end{assumption}

\begin{assumption} \label{ass:2}
 The function $g(x,y)$ is nonconvex and satisfies
\begin{align}
 {\color{blue}{\varrho\big(\nabla^2_{yy}g\big(x,y^*(x)\big)\big) \in [\mu, L_g]}},
\end{align}
where $y^*(x) \in \arg\min_y g(x,y)$, and $\varrho(\cdot)$ denotes the eigenvalue (or singular-value) function and $L_g\geq\mu>0$.
\end{assumption}

\begin{assumption} \label{ass:3}
The functions $f(x,y)$ and $g(x,y)$ satisfy
\begin{itemize}
 \item[(i)] For the minimizers $y^*(x) \in \arg\min_y g(x,y)$, we have
 \begin{align}
 {\color{blue}{\|\nabla_yf(x,y^*(x))\|\leq C_{fy}, \ \|\nabla^2_{xy}g(x,y^*(x))\|\leq C_{gxy}}};
 \end{align}
 \item[(ii)] The partial derivatives $\nabla_x f(x,y)$ and $\nabla_y f(x,y)$ are $L_f$-Lipschitz continuous;
 \item[(iii)] The partial derivatives $\nabla_xg(x,y)$ and $\nabla_yg(x,y)$ are $L_g$-Lipschitz continuous.
\end{itemize}
\end{assumption}

 Assumption~\ref{ass:1} is commonly used in bilevel optimization without the lower-level strongly convexity~\citep{liu2022bome,shen2023penalty,huang2023momentum}.
 Assumption~\ref{ass:2} is a milder condition in bilevel optimization without the lower-level strongly-convexity,
 which imposes the non-singularity of $\nabla^2_{yy}g(x,y)$ only at the minimizers $y^*(x) \in \arg\min_y g(x,y)$, as in~\citep{huang2023momentum}. \textbf{Note that} since $y^*(x)\in \arg\min_y g(x,y)$, we can not have negative eigenvalues at the minimizer $y^*(x)$, so Assumption~2 assumes that $\varrho\big(\nabla^2_{yy}g\big(x,y^*(x)\big)\big) \in [\mu, L_g]$ instead of $\varrho\big(\nabla^2_{yy}g\big(x,y^*(x)\big)\big) \in [-L_g,-\mu] \cup [\mu, L_g]$.
Since $\nabla^2_{yy} g(x,y)$
is a symmetric matrix, its singular values are the absolute value of eigenvalues.
Hence, we also can use $\varrho(\cdot)$ to denote the singular-value function.

\textbf{Note that} the GALET~\citep{xiao2023generalized} method simultaneously uses the PL condition, its Assumption 2 (i.e., let $\sigma_g = \inf_{x,y}\{\sigma_{\min}^{+}(\nabla^2_{yy} g(x,y))\} >0$ for all $(x,y)$) and its Assumption 1 (i.e., $\nabla^2_{yy} g(x,y)$ is Lipschitz continuous). Clearly, when Hessian matrix $\nabla^2_{yy} g(x,y)$ is singular, its Assumption 1 and Assumption 2 imply that the lower bound of the non-zero singular values $\sigma_g$ is close to zero (i.e., $\sigma_g\rightarrow 0$), under this case, the convergence results of the GALET are \textbf{meaningless}, e.g., the constant $L_w = \frac{\ell_{f,1}}{\sigma_g}+\frac{\sqrt{2}\ell_{g,2}\ell_{f,0}}{\sigma_g^2}\rightarrow + \infty$ used in its Lemmas 6 and 9, and $L_F = \ell_{f,0}(\ell_{f,1}+\ell_{g,2})/\sigma_g \rightarrow + \infty$ used in its Lemma 12. Under the other case, the PL condition, Lipschitz continuous of Hessian and its Assumption 2 (the singular values of Hessian is bounded away from 0, i.e., $\sigma_g>0$) imply that the GALET uses the strongly convex assumption. \textbf{Note that} although the singular values of Hessian $\nabla^2_{yy} g(x,y)$ \textbf{exclude} zero, i.e, the eigenvalues of Hessian $\nabla^2_{yy} g(x,y)$ may be in $[-\ell_{g,2},-\sigma_g]\bigcup [\sigma_g,\ell_{g,2}]$, we
cannot have negative eigenvalues at the minimizer $y^*(x)$. Meanwhile, since Hessian is Lipschitz continuous, its all eigenvalues are in $[\sigma_g,\ell_{g,2}]$. Thus, the PL condition, Lipschitz continuous of Hessian and its Assumption 2 (the singular values of Hessian is bounded away from 0, i.e., $\sigma_g>0$) imply that the GALET assumes the strongly convex.

For example, considering a PL function: $g(x,y) = y^2+\sin^2(y)+x$ on variable $y$, where both $x$ and $y$ are scalar variables, we have $\nabla^2_{yy} g(x,y)=2+2\cos^2(y)-2\sin^2(y)$ and
$0=y^*=\arg\min_y g(x,y)$, and then we get $\nabla^2_{yy} g(x,y^*)=\nabla^2_{yy} g(x,0) =4>0$ and $\nabla^2_{yy} g(x,\pi/2)=0$. Clearly, Assumption~2 of the GALET~\citep{xiao2023generalized} can not conform to this PL function, since $\nabla^2_{yy} g(x,\pi/2)=0$. Clearly, our Assumption~\ref{ass:2} is milder than the Assumption~2 of the GALET. Meanwhile, our Assumption~\ref{ass:2} is more reasonable in practice.

Assumption~\ref{ass:3} is commonly appeared in bilevel optimization methods \citep{ghadimi2018approximation,ji2021bilevel,liu2022bome}.
 \textbf{Note that} our Assumption~3(i) assumes that
 $\|\nabla_yf(x,y)\|$ and $\|\nabla^2_{xy}g(x,y)\|$ are bounded \textbf{only at
 the minimizer} $y^*(x) \in \arg\min_y g(x,y)$ as in \citep{huang2023momentum}, while \citep{ghadimi2018approximation,ji2021bilevel} assume that $\|\nabla_yf(x,y)\|$ and $\|\nabla^2_{xy}g(x,y)\|$ are bounded \textbf{at
 any} $y \in \mathbb{R}^p$ and $x\in \mathbb{R}^d$.
 Meanwhile, the BOME~\cite{liu2022bome} uses the stricter assumption that
  $\|\nabla f(x,y)\|$, $\|\nabla g(x,y)\|$, $|f(x,y)|$ and $|g(x,y)|$
 are bounded for any $(x,y)$ in its Assumption~3.
 More recently, \cite{liu2023averaged} claim that
 the proposed sl-BAMM~\citep{liu2023averaged} method does not require the bounded $\nabla_y F(x,y)$,
however, it needs to $\|x\|+\|y\|\leq D$ for all $x\in \mathbb{R}^n, y\in \mathbb{R}^m$ and the bounded $\nabla_y F(0,0)$. For solving the unconstrained bilevel optimization problems, clearly, the condition $\|x\|+\|y\|\leq D$ for all $x\in \mathbb{R}^n, y\in\mathbb{R}^m$ may be not satisfied.

For example, considering a bilevel problem $\min_{x\in [0,1],y\in y^*(x)} \Big\{f(x,y)=xy^2+y-x-1+x^3\Big\}, \ \mbox{s.t.} y^*(x) \equiv\min_{y\in \mathbb{R}} \Big\{g(x,y)=y^2+\sin^2(y)+x\Big\}$, which can be rewritten as $\min_{x\in \mathbb{R},y\in y^*(x)} \Big\{f(x,y)+ h(x)\Big\}, \ \mbox{s.t.} y^*(x) \equiv\min_{y\in \mathbb{R}} \Big\{g(x,y)=y^2+\sin^2(y)+x\Big\}$ where $h(x)=0$ when $x\in [0,1]$ otherwise $h(x)=+\infty$,
clearly, we can obtain $y^*(x)=0$, $\nabla^2_{yy}g(x,y)=2+2\cos^2(y)-2\sin^2(y)$, $\nabla^2_{xy}g(x,y)=0$, $\nabla_y f(x,y)=2xy+1$, $\nabla^2_{yy}g(x,y^*(x))=4$, $\nabla^2_{xy}g(x,y^*(x))=0$ and $\nabla_y f(x,y^*(x))=1$. Since $\nabla^2_{yy}g(x,y^*(x))=4$, $\nabla^2_{xy}g(x,y^*(x))=0$ and $\nabla_y f(x,y^*(x))=1$, our Assumption~\ref{ass:2} and Assumption~\ref{ass:3}(i) holds. Meanwhile  $\nabla^2_{yy}g(x,y)=2+2\cos^2(y)-2\sin^2(y)$ for any $y\in \mathbb{R}$ may be zero or negative such as $\nabla^2_{yy} g(x,\pi/2)=0$, and $\nabla_y f(x,y)=2xy+1$ may be not bounded for any $y\in \mathbb{R}$. Thus, the assumptions used in~\citep{liu2022bome,xiao2023generalized} may be not satisfied.

\begin{assumption} \label{ass:4}
 The partial derivatives $\nabla^2_{xy}g(x,y)$ and $\nabla^2_{yy}g(x,y)$ are $L_{gxy}$-Lipschitz and $L_{gyy}$-Lipschitz, e.g.,
 for all $x,x_1,x_2 \in \mathbb{R}^d$ and $y,y_1,y_2 \in \mathbb{R}^p$
 \begin{align}
   \|\nabla^2_{xy} g(x_1,y)-\nabla^2_{xy} g(x_2,y)\| \leq L_{gxy}\|x_1-x_2\|,  \ \|\nabla^2_{xy} g(x,y_1)-\nabla^2_{xy} g(x,y_2)\| \leq L_{gxy}\|y_1-y_2\|. \nonumber
 \end{align}
 \end{assumption}

\begin{assumption} \label{ass:5}
 The function $F(x)=f(x,y^*(x))$ is bounded below in $x\in \mathbb{R}^d$, \emph{i.e.,} $F^* = \inf_{x\in \mathbb{R}^d}F(x) > -\infty$.
\end{assumption}

Assumption~\ref{ass:4} is also commonly used in bilevel optimization methods \citep{ghadimi2018approximation,ji2021bilevel}.
Assumption~\ref{ass:5} ensures the feasibility of the bilevel Problems~(\ref{eq:1}) and~(\ref{eq:2}).

\subsection{Useful Lemmas}
In this subsection, based on the above assumptions, we give some useful lemmas.

\begin{lemma} \label{lem:1}
(\cite{huang2023momentum})
Under the above Assumption \ref{ass:2}, we have, for any $x\in \mathbb{R}^d$,
\begin{align}
 \nabla F(x) & =\nabla_x f(x,y^*(x)) - \nabla^2_{xy} g(x,y^*(x))\Big[\nabla^2_{yy}g(x,y^*(x))\Big]^{-1}\nabla_y f(x,y^*(x)).
\end{align}
\end{lemma}

From the above Lemma~\ref{lem:1}, we can get the form of hyper-gradient $\nabla F(x)$ is the same as the hyper-gradient in~(\ref{eq:g1}). Since the Hessian matrix $\nabla^2_{yy}g(x,y)$ for all $(x,y)$ may be singular and both $\nabla^2_{xy}g(x,y)$ and $\nabla_yf(x,y)$ are not bounded for all $(x,y)$, as in \cite{huang2023momentum},
we define a useful hyper-gradient estimator:
\begin{align}
\hat{\nabla} f(x,y)=\nabla_xf(x,y) - {\color{blue}{\hat{\Pi}_{C_{gxy}}\big[\nabla^2_{xy}g(x,y)\big]}} \big({\color{blue}{\mathcal{S}_{[\mu,L_g]}\big[\nabla^2_{yy}g(x,y)\big]}}\big)^{-1}
{\color{blue}{\Pi_{C_{fy}}\big[\nabla_yf(x,y)\big]}},
\end{align}
which replaces the standard hyper-gradient estimator $\breve{\nabla} f(x,y)$ used in ~\cite{ghadimi2018approximation,ji2021bilevel} for the strongly-convex lower-level problems:
\begin{align}
\breve{\nabla} f(x,y)=\nabla_xf(x,y) - \nabla^2_{xy}g(x,y) \big(\nabla^2_{yy}g(x,y)\big)^{-1}\nabla_yf(x,y).
\end{align}

\begin{lemma} \label{lem:2}
(\cite{huang2023momentum})
Under the above Assumptions \ref{ass:1}-\ref{ass:4}, the functions (or mappings) $F(x)=f(x,y^*(x))$, $G(x)=g(x,y^*(x))$ and $y^*(x)\in \arg\min_{y\in \mathbb{R}^p}g(x,y)$ satisfy, for all $x_1,x_2\in \mathbb{R}^d$,
\begin{align}
 & \|y^*(x_1)-y^*(x_2)\| \leq \kappa\|x_1-x_2\|, \quad \|\nabla y^*(x_1) - \nabla y^*(x_2)\| \leq L_y\|x_1-x_2\| \nonumber \\
 & \|\nabla F(x_1) - \nabla F(x_2)\|\leq L_F\|x_1-x_2\|, \quad \|\nabla G(x_1) - \nabla G(x_2)\|\leq L_G\|x_1-x_2\| \nonumber
\end{align}
where $\kappa=C_{gxy}/\mu$, $L_y=\big( \frac{C_{gxy}L_{gyy}}{\mu^2} +  \frac{L_{gxy}}{\mu} \big) (1+ \frac{C_{gxy}}{\mu})$,
$L_F=\Big(L_f + L_f\kappa + C_{fy}\big( \frac{C_{gxy}L_{gyy}}{\mu^2} +  \frac{L_{gxy}}{\mu} \big)\Big)(1+\kappa)$ and $L_G=\Big(L_g + L_g\kappa + C_{gy}\big( \frac{C_{gxy}L_{gyy}}{\mu^2} +  \frac{L_{gxy}}{\mu} \big)\Big)(1+\kappa)$.
\end{lemma}

\begin{lemma} \label{lem:3}
(\cite{huang2023momentum})
Let $\hat{\nabla} f(x,y)=\nabla_xf(x,y) - {\color{blue}{\hat{\Pi}_{C_{gxy}}\big[\nabla^2_{xy}g(x,y)\big] }} \big({\color{blue}{\mathcal{S}_{[\mu,L_g]}\big[\nabla^2_{yy}g(x,y)\big]}}\big)^{-1}{\color{blue}{\Pi_{C_{fy}}\big[\nabla_yf(x,y)\big]}}$ and
$\nabla F(x_t) = \nabla f(x_t,y^*(x_t))$,
 we have
 \begin{align}
 \|\hat{\nabla} f(x,y)-\nabla F(x)\|^2 \leq \hat{L}^2\|y^*(x)-y\|^2\leq \frac{2\hat{L}^2}{\mu}\big(g(x,y)-\min_y g(x,y)\big),
\end{align}
where $\hat{L}^2 = 4\big(L^2_f+ \frac{L^2_{gxy}C^2_{fy}}{\mu^2} + \frac{L^2_{gyy} C^2_{gxy}C^2_{fy}}{\mu^4} +
 \frac{L^2_fC^2_{gxy}}{\mu^2}\big)$.
\end{lemma}

\section{ Adaptive Mirror Descent Bilevel Algorithms }
In the section, we study a class of nonconvex bilevel optimization problems, and
propose an efficient adaptive projection-aid gradient (i.e., AdaPAG) method based on mirror descent to solve the deterministic Problem~(\ref{eq:1}).
To solve the stochastic Problem~(\ref{eq:2}), further we propose an efficient adaptive stochastic projection-aid gradient (i.e., AdaVSPAG) methods based on mirror descent and variance-reduced techniques.

We first review the mirror distance and mirror descent~\citep{censor1981iterative,censor1992proximal,chen1993convergence,beck2003mirror}. Given a differential and (strongly) convex Bregman function $\psi(\cdot)$, the mirror distance can be defined as
\begin{align}
 D(x,x_0) = \psi(x) - \psi(x_0) -\langle\nabla \psi(x_0), x-x_0\rangle, \ \forall x, x_0 \in \mathbb{R}^d
\end{align}
When solving the minimization problem $\min_{x\in \mathbb{R}}f(x)$, at $t$-th iteration, the mirror descent can be represented as
\begin{align}
 x_{t+1} = \arg\min_{x\in \mathbb{R}^d} \Big\{ \langle \nabla f(x_t), x-x_t\rangle + \frac{1}{\gamma}D(x,x_t)\Big\},
\end{align}
where $\gamma>0$ denotes the learning rate. In the paper, we consider a generalized mirror descent, also called adaptive mirror descent, can be defined as
\begin{align}
 x_{t+1} = \arg\min_{x\in \mathbb{R}^d} \Big\{ \langle \nabla f(x_t), x-x_t\rangle + \frac{1}{\gamma}D_t(x,x_t)\Big\},
\end{align}
where $\gamma>0$ denotes the learning rate, and $D_t(x,x_t)$ is an adaptive mirror distance,
\begin{align}
 D_t(x,x_0) = \psi_t(x) - \psi_t(x_0) -\langle\nabla \psi_t(x_0), x-x_0\rangle, \ \forall x, x_0 \in \mathbb{R}^d
\end{align}
where $\psi_t(\cdot)$ is an adaptive (dynamic) Bregman function, which may be depend on
the value functions and gradients $\{f(x_i),\nabla f(x_i)\}_{i=1}^t$. For example,
in the adaptive gradient methods, we can set adaptive Bregman function $\psi_t(x)=\frac{1}{2}x^TA_tx$ with
$A_t\succ 0$, which is a matrix form of adaptive learning rate. In the popular Adam algorithm~\citep{kingma2014adam}, $A_t$ is a diagonal matrix, i.e., $A_t(i,i)=\sqrt{\hat{v}_t(i)+\epsilon}$ with $\epsilon>0$ for all $i\in [d]$, and $A(i,j)=0$ for all $i\neq j$.

\begin{algorithm}[tb]
\caption{ AdaPAG Algorithm }
\label{alg:1}
\begin{algorithmic}[1] 
\STATE {\bfseries Input:} $T$, parameters $\gamma>0$, $\lambda>0$
and initial input $x_1 \in \mathbb{R}^d$ and $y_1 \in \mathbb{R}^p$; \\
\FOR{$t = 1, 2, \ldots, T$}
\STATE $v_t = \nabla_y g(x_t,y_t)$, $u_t = \nabla_x f(x_t,y_t)$, $h_t =\Pi_{C_{fy}}\big[\nabla_yf(x_t,y_t)\big]$, $G_t = \hat{\Pi}_{C_{gxy}}\big[\nabla^2_{xy}g(x_t,y_t)\big]$;
\STATE  {\color{blue}{ $H_t = \mathcal{S}_{[\mu,L_g]}\big[\nabla_{yy}^2 g(x_t,y_t)\big] = U_t\Theta_tU_t^T$ }}, where $\theta_{t,i}\in [\mu,L_g]$ for all $i=1,\cdots,p$;
\STATE $w_t=u_t - G_t(H_t)^{-1}h_t= \nabla_x f(x_t,y_t) - G_t\Big(\sum_{i=1}^p \big(U^T_{t,i}h_t \big)/\theta_{t,i}U_{t,i}\Big)$;
\STATE Generate the adaptive matrices $A_t \in \mathbb{R}^{d \times d}$ and $B_t \in \mathbb{R}^{p \times p}$;\\
\textcolor{blue}{One example of $A_t$ and $B_t$ by using update rule ($a_0 = 0$, $b_0 = 0$, $ 0 < \tau< 1$, $\rho>0$) } \\
\textcolor{blue}{ Compute $ a_t = \tau a_{t-1} + (1 - \tau)v_t^2$, $A_t = \mbox{diag}(\sqrt{a_t} + \rho)$}; \\
\textcolor{blue}{ Compute $ b_t = \tau b_{t-1} + (1 - \tau)w_t^2$, $B_t = \mbox{diag}(\sqrt{b_t} + \rho)$}; \\
\STATE $x_{t+1} = \arg\min_{x\in \mathbb{R}^d}\big\{ \langle w_t, x\rangle + \frac{1}{2\gamma}(x-x_t)^TA_t(x-x_t) + h(x) \big\}$;
\STATE $y_{t+1} = \arg\min_{y\in \mathbb{R}^p}\big\{ \langle v_t, y\rangle + \frac{1}{2\lambda}(y-y_t)^TB_t(y-y_t) \big\}$;
\ENDFOR
\STATE {\bfseries Output:} Chosen uniformly random from $\{x_t\}_{t=1}^{T}$.
\end{algorithmic}
\end{algorithm}

\subsection{Adaptive Projection-Aid Gradient Algorithm}
In this subsection, we propose an efficient adaptive projection-aid gradient (i.e., AdaPAG) method based on mirror descent to solve the problem~(\ref{eq:1}). Algorithm \ref{alg:1} provides a procedure framework of our AdaPAG algorithm.

At the line 5 of Algorithm~\ref{alg:1}, we estimate the gradient $\nabla F(x_t)$ based on the above Lemma~\ref{lem:1},
\begin{align}
w_t = \hat{\nabla}f(x_t,y_t) =\nabla_xf(x_t,y_t) + \hat{\Pi}_{C_{gxy}}\big[\nabla^2_{xy}g(x_t,y_t)\big] \big(\mathcal{S}_{[\mu,L_g]}\big[\nabla^2_{yy}g(x_t,y_t)\big]\big)^{-1}\Pi_{C_{fy}}\big[\nabla_yf(x_t,y_t)\big], \nonumber
\end{align}
where the projection operator $\mathcal{S}_{[\mu,L_g]}$ is implemented by using SVD and thresholding the singular values. Then we have
\begin{align}
 H_t = \mathcal{S}_{[\mu,L_g]}\big[\nabla_{yy}^2 g(x_t,y_t)\big] = U_t\Theta_tU_t^T,
\end{align}
where $U_t$ is a real orthogonal matrix, i.e, $U_t^TU_t=U_tU_t^T=I_p$, and $\Theta_t=\mbox{diag}(\theta_t) \in \mathbb{R}^{p\times p}$ is a diagonal matrix with $\theta_{t,i} \in [\mu,L_g]$ for all $i=1,2,\cdots, p$. Clearly, we can easily obtain
\begin{align}
 H_t^{-1} = U_t\Theta_t^{-1}U_t^T, \quad \Theta_t^{-1} = \mbox{diag}(1/\theta_t),
\end{align}
where $1/\theta_t = (1/\theta_{t,1},\cdots,1/\theta_{t,p})$.
Thus, we can get
\begin{align} \label{eq:3}
 w_t=u_t - G_t(H_t)^{-1}h_t & = \nabla_x f(x_t,y_t) - G_tU_t\Theta_t^{-1}U_t^Th_t \nonumber \\
 & = \nabla_x f(x_t,y_t) - G_t \Big(\sum_{i=1}^p \big(U^T_{t,i}h_t \big)/\theta_{t,i}U_{t,i}\Big),
\end{align}
where $U_{t,i}$ is the $i$-th column of matrix $U_t = [U_{t,1},\cdots, U_{t,p}]$.
Based on the above equality~(\ref{eq:3}), calculations of the gradient $w_t$ mainly depends on the above SVD.
Recently, the existing many methods~\citep{tzeng2013split} can compute SVD of the large-scale matrices.
Thus, our algorithm can also solve the large-scale bilevel problems in high-dimension setting.

Although our Algorithm~\ref{alg:1} requires to compute Hessian matrix, its inverse and
its projecting matrix over spectral set, our algorithm can avoid computing inverse of Hessian matrix. In fact, many bilevel optimization methods~\citep{ghadimi2018approximation,ji2021bilevel,xiao2023generalized} still compute Hessian matrix $\nabla^2_{yy}g(x,y)$.
\textbf{Note that} the GALET~\citep{xiao2023generalized} method not only calculates Hessian matrices, but also computes the product of Hessian matrices multiple times when updating variable $x$ at each iteration.

In our AdaPAG algorithm, we use the adaptive Bregman functions $\psi_t(x)=\frac{1}{2}x^TA_tx$ and
$\phi_t(x)=\frac{1}{2}x^TB_tx$ to generate adaptive mirror distances in updating the variables $x$ and $y$, respectively.
For example, $\psi_t(x)=\frac{1}{2}x^TA_tx$, the adaptive mirror distance can be defined as
\begin{align}
  D_t(x,x_t) = \psi_t(x) - \psi_t(x_t) -\langle\nabla \psi_t(x_t), x-x_t\rangle
  = \frac{1}{2}(x-x_t)^TA_t(x-x_t) .
\end{align}
At the line 7 of Algorithm~\ref{alg:1}, we have
\begin{align}
x_{t+1} & = \arg\min_{x\in \mathbb{R}^d}\Big\{ \langle w_t, x-x_t\rangle
+ \frac{1}{2\gamma}(x-x_t)^TA_t(x-x_t) + h(x) \Big\} \nonumber \\
& = \arg\min_{x\in \mathbb{R}^d}\Big\{ \frac{1}{2\gamma} \|x-(x_t- A_t^{-1}w_t)\|^2_{A_t} + h(x) \Big\} = \mathbb{P}_{h(\cdot)}^\gamma(x_t - A_t^{-1}w_t),
\end{align}
where $\|z\|^2_{A_t} = z^TA_tz$ for all $z\in \mathbb{R}^d$.

At the line 8 of Algorithm~\ref{alg:1}, we have
\begin{align}
y_{t+1} = y_t - \lambda B_t^{-1}v_t = \arg\min_{y\in \mathbb{R}^p}\Big\{ \langle v_t, y-y_t\rangle
+ \frac{1}{2\lambda}(x-x_t)^TB_t(x-x_t) \Big\},
\end{align}
where $\lambda>0$ is a hyper-parameter, and $B_t$ is an adaptive matrix. Note that these
adaptive matrices $A_t$ and $B_t$ generally the diagonal matrices. So the inverse of matrices $A_t$ and $B_t$
can not induce high computation.
When using the adaptive matrix $B_t$ generated from Algorithm~\ref{alg:1}, we have
\begin{align}
y_{t+1,i} = y_{t,i} - \lambda/(b_{t,i}+\rho) v_{t,i}, \quad \forall i\in [p].
\end{align}
where $y_{t,i}$ denotes the $i$-th element of vector $y_t\in \mathbb{R}^p$.
Meanwhile, we can also choose some other adaptive matrices $A_t$ and $B_t$ generated as the Adagrad algorithm~\citep{duchi2011adaptive}
such as
\begin{align}
 A_t = \frac{1}{\sqrt{\sum_{k=1}^t||w_k||^2}}I_d, \quad B_t = \frac{1}{\sqrt{\sum_{k=1}^t||v_k||^2}}I_p.
\end{align}

 \begin{algorithm}[tb]
\caption{ AdaVSPAG Algorithm }
\label{alg:2}
\begin{algorithmic}[1] 
\STATE {\bfseries Input:} $T$, stepsizes $\gamma>0$, $\lambda>0$, mini-batch sizes $q$, $n$, $b$; \\
\STATE {\bfseries initialize:}  $x_0 \in \mathbb{R}^d$ and $y_0 \in \mathbb{R}^{p}$;  \\
\FOR{$t = 0, 1, \cdots, T-1$}
\IF {$\mod(t,q)=0$}
\STATE Draw $2n$ independent samples $\{\xi^i_t\}_{i=1}^n$ and $\{\zeta^i_t\}_{i=1}^n$;
\STATE $v_t = \frac{1}{n}\sum_{i=1}^n\nabla_y g (x_t,y_t;\zeta^i_t)$, $G_t = \hat{\Pi}_{C_{gxy}}\big[\frac{1}{n}\sum_{i=1}^n \nabla_{xy}^2 g(x_t,y_t;\zeta^i_t)\big]= W_t\Sigma_tV_t^T$; \\
\STATE {\color{blue}{ $H_t =\mathcal{S}_{[\mu,L_g]}\big[\frac{1}{n}\sum_{i=1}^n\nabla_{yy}^2 g(x_t,y_t;\zeta^i_t)\big] = U_t\Theta_tU_t^T$ }}; \\
\STATE $u_t = \frac{1}{n}\sum_{i=1}^n\nabla_x f(x_t,y_t;\xi^i_t)$, $h_t = \Pi_{C_{fy}}\big[ \frac{1}{n}\sum_{i=1}^n\nabla_y f(x_t,y_t;\xi^i_t) \big]$; \\
\STATE $w_t = u_t - G_t(H_t)^{-1}h_t = u_t - W_t\Sigma_tV_t^T\Big(\sum_{i=1}^p \big(U^T_{t,i}h_t \big)/\theta_{t,i}U_{t,i}\Big)$; \\
\ELSE
\STATE Draw $2b$ independent samples $\{\xi^i_t\}_{i=1}^b$ and $\{\zeta^i_t\}_{i=1}^b$;
\STATE $v_t = \frac{1}{b}\sum_{i=1}^b\nabla_y g (x_t,y_t;\zeta^i_t) - \frac{1}{b}\sum_{i=1}^b\nabla_y g (x_{t-1},y_{t-1};\zeta^i_t) + v_{t-1}$; \\
\STATE $G_t = \hat{\Pi}_{C_{gxy}}\big[\frac{1}{b}\sum_{i=1}^b \nabla_{xy}^2 g(x_t,y_t;\zeta^i_t)
 - \frac{1}{b}\sum_{i=1}^b\nabla_{xy}^2 g (x_{t-1},y_{t-1};\zeta^i_t) + G_{t-1}\big]= W_t\Sigma_tV_t^T$; \\
\STATE {\color{blue}{$H_t =\mathcal{S}_{[\mu,L_g]}\big[\frac{1}{b}\sum_{i=1}^b\nabla_{yy}^2 g(x_t,y_t;\zeta^i_t)- \frac{1}{b}\sum_{i=1}^b\nabla_{yy}^2 g (x_{t-1},y_{t-1};\zeta^i_t) + H_{t-1}\big]= U_t\Theta_tU_t^T$ }}; \\
\STATE $u_t = \frac{1}{b}\sum_{i=1}^b\nabla_x f(x_t,y_t;\xi^i_t) - \frac{1}{b}\sum_{i=1}^b\nabla_x f(x_{t-1},y_{t-1};\xi^i_t) + u_{t-1}$; \\
\STATE $h_t = \Pi_{C_{fy}}\big[ \frac{1}{b}\sum_{i=1}^b\nabla_y f(x_t,y_t;\xi^i_t) - \frac{1}{b}\sum_{i=1}^b\nabla_y f(x_{t-1},y_{t-1};\xi^i_t) + h_{t-1}\big]$; \\
\STATE $w_t = u_t - G_t(H_t)^{-1}h_t= u_t - W_t\Sigma_tV_t^T\Big(\sum_{i=1}^p \big(U^T_{t,i}h_t \big)/\theta_{t,i}U_{t,i}\Big)$; \\
\ENDIF
\STATE Generate the adaptive matrices $A_t \in \mathbb{R}^{d \times d}$ and $B_t \in \mathbb{R}^{p \times p}$;\\
\textcolor{blue}{One example of $A_t$ and $B_t$ by using update rule ($a_0 = 0$, $b_0 = 0$, $ 0 < \tau< 1$, $\rho>0$) } \\
\textcolor{blue}{ Compute $ a_t = \tau a_{t-1} + (1 - \tau)v_t^2$, $A_t = \mbox{diag}(\sqrt{a_t} + \rho)$}; \\
\textcolor{blue}{ Compute $ b_t = \tau b_{t-1} + (1 - \tau)w_t^2$, $B_t = \mbox{diag}(\sqrt{b_t} + \rho)$}; \\
\STATE $x_{t+1} = \arg\min_{x\in \mathbb{R}^d}\big\{ \langle w_t, x\rangle
+ \frac{1}{2\gamma}(x-x_t)^TA_t(x-x_t) + h(x)\big\}$;
\STATE $y_{t+1} = \arg\min_{y\in \mathbb{R}^p}\big\{ \langle v_t, y\rangle
+ \frac{1}{2\lambda}(y-y_t)^TB_t(y-y_t) \big\}$;
\ENDFOR
\STATE {\bfseries Output:} Uniformly and randomly choose from $\{x_t\}_{t=1}^{T}$.
\end{algorithmic}
\end{algorithm}

\subsection{ Adaptive Stochastic Projection-Aid Gradient Algorithm }
In this subsection, we present an efficient adaptive stochastic projection-aid gradient (i.e., AdaVSPAG) methods based on mirror descent and variance-reduced techniques to solve the problem~(\ref{eq:2}).
Algorithm \ref{alg:2} shows a procedure framework of our AdaVSPAG algorithm.

In our AdaVSPAG algorithm, we use the variance reduced technique of
SARAH/SPIDER~\citep{nguyen2017sarah,fang2018spider,wang2019spiderboost} to estimate stochastic gradients:
\begin{align}
& v_t = \frac{1}{b}\sum_{i=1}^b\big(\nabla_y g (x_t,y_t;\zeta^{1,i}_t) - \nabla_y g (x_{t-1},y_{t-1};\zeta^{1,i}_t) \big) +v_{t-1}\nonumber \\
& u_t = \frac{1}{b}\sum_{i=1}^b\nabla_x f(x_t,y_t;\xi^i_t) - \frac{1}{b}\sum_{i=1}^b\nabla_x f(x_{t-1},y_{t-1};\xi^i_t) + u_{t-1}, \nonumber
\end{align}
and stochastic projected gradients:
\begin{align}
&G_t = \hat{\Pi}_{C_{gxy}}\big[\frac{1}{b}\sum_{i=1}^b \nabla_{xy}^2 g(x_t,y_t;\zeta^i_t)
 - \frac{1}{b}\sum_{i=1}^b\nabla_{xy}^2 g (x_{t-1},y_{t-1};\zeta^i_t) + G_{t-1}\big] \nonumber \\
& h_t = \Pi_{C_{fy}}\big[ \frac{1}{b}\sum_{i=1}^b\nabla_y f(x_t,y_t;\xi^i_t) - \frac{1}{b}\sum_{i=1}^b\nabla_y f(x_{t-1},y_{t-1};\xi^i_t) + h_{t-1}\big] \nonumber \\
& H_t =\mathcal{S}_{[\mu,L_g]}\big[\frac{1}{b}\sum_{i=1}^b\nabla_{yy}^2 g(x_t,y_t;\zeta^i_t)- \frac{1}{b}\sum_{i=1}^b\nabla_{yy}^2 g (x_{t-1},y_{t-1};\zeta^i_t) + H_{t-1}\big].\nonumber
\end{align}

In our AdaVSPAG algorithm, in updating variables $x$,
we also use the following adaptive mirror descent:
\begin{align}
x_{t+1}  = \arg\min_{x\in \mathbb{R}^d}\Big\{ \langle w_t, x\rangle
+ \frac{1}{2\gamma}(x-x_t)^TA_t(x-x_t) + h(x) \Big\}, \nonumber
\end{align}
where $\gamma>0$ is a learning rate.
It is similar for updating variable $y$.

\section{Convergence Analysis}
In the section, we study the convergence properties of our algorithms under some mild assumptions.
We first define a useful gradient mapping
$\mathcal{G}(x_t,\nabla F(x_t),\gamma)=\frac{1}{\gamma}(x_t-x_{t+1})$, where the sequence $\{x_t\}_{t\geq1}$
is generated from
\begin{align}
x_{t+1} & = \arg\min_{x\in \mathbb{R}^d}\Big\{ \langle \nabla F(x_t), x\rangle
+ \frac{1}{2\gamma}(x-x_t)^TA_t(x-x_t) + h(x)\Big\} \\
& = \arg\min_{x\in \mathbb{R}^d}\Big\{ \frac{1}{2\gamma} \|x-(x_t- A_t^{-1}\nabla F(x_t))\|^2_{A_t} + h(x)\Big\} = \mathbb{P}_{h(\cdot)}^\gamma(x_t - A_t^{-1}\nabla F(x_t)),
\end{align}
where $F(x_t)=f(x,y^*(x))$ with $y^*(x)\in \arg\min_y g(x,y)$ and $\|z\|^2_{A_t} = z^TA_tz$ for all $z\in \mathbb{R}^d$.

\begin{assumption} \label{ass:6}
In our algorithms, the adaptive matrices $A_t$ and $B_t$ for all $t\geq 1$ satisfy $A_t\succeq \rho I_d \succ 0$ and $\rho_u I_{p} \succeq B_t \succeq \rho_l I_{p} \succ 0$ for any $t\geq 1$, where $\rho >0,\rho_u \geq \rho_l >0$ are appropriate positive numbers.
\end{assumption}

Assumption \ref{ass:6} ensures that the adaptive matrices $A_t$ and $B_t$ for all $t\geq 1$ are positive definite as in \citep{huang2021super}. In fact, the adaptive Bregman functions $\psi_t(x)=\frac{1}{2}x^TA_tx$ and
$\phi_t(x)=\frac{1}{2}x^TB_tx$ are smooth and strongly convex imply that Assumption \ref{ass:6} holds. Thus, the above Assumption \ref{ass:6} is mild.
\textbf{Note that} the adaptive bilevel methods in \cite{huang2021biadam} only consider using the global adaptive learning rates (i.e., $B_t=b_t I_p \ (b_t>0)$) to update the dual variable $y$.

\subsection{Convergence Analysis of AdaPAG Algorithm}
In this subsection,
we provide the convergence properties of our AdaPAG algorithm under some mild assumptions.

\subsubsection{Convergence Properties on Unimodal $g(x,y)$}
We first present the convergence properties of our AdaPAG method when $g(x,\cdot)$ satisfies the global PL condition for all $x$, i.e., it has \textbf{a unique minimizer} $y^*(x)=\arg\min_y f(x,y)$. We first give a useful lemma.

\begin{lemma} \label{lem:4}
Suppose the sequence $\{x_t,y_t\}_{t=1}^T$ be generated from Algorithms~\ref{alg:1} or \ref{alg:2}.
Under the above Assumptions~\ref{ass:1}-\ref{ass:3} and \ref{ass:6}, given $\gamma\leq \min\Big\{\frac{\lambda\mu\rho_l}{16L_G\rho_u}, \frac{\mu\rho_l}{16L^2_g\rho_u} \Big\}$, $0<\rho_l\leq 1$ and $0<\lambda \leq \frac{1}{2L_g\rho_u}$, we have
\begin{align}
g(x_{t+1},y_{t+1}) - G(x_{t+1})
& \leq (1-\frac{\lambda\mu}{2\rho_u}) \big(g(x_t,y_t) -G(x_t)\big) + \frac{\rho_l}{8\gamma}\|x_{t+1}-x_t\|^2  -\frac{1}{4\lambda\rho_u}\|y_{t+1}-y_t\|^2 \nonumber \\
& \qquad + \frac{\lambda}{\rho_l}\|\nabla_y g(x_t,y_t)-v_t\|^2,
\end{align}
where $G(x_t)=g(x_t,y^*(x_t))$ with $y^*(x_t) \in \arg\min_{y}g(x_t,y)$ for all $t\geq 1$.
\end{lemma}
The above Lemma~\ref{lem:4} shows the properties of the residuals $g(x_t,y_t)-G(x_t)\geq 0$ for all $t\geq 1$. In fact, this lemma can also be used in the following Algorithms~\ref{alg:2}.

\begin{theorem}  \label{th:1}
 Assume the sequence $\{x_t,y_t\}_{t=1}^T$ be generated from our Algorithm \ref{alg:1}. Under the above Assumptions~\ref{ass:1}-\ref{ass:6}, let $0< \gamma \leq \min\big(\frac{3\rho}{4L_F},\frac{\rho\lambda\mu^2}{8\rho_u\hat{L}^2},\frac{\lambda\mu\rho_l}{16L_G\rho_u}, \frac{\mu\rho_l}{16L^2_g\rho_u} \big)$, $0<\lambda\leq \frac{1}{2L_g\rho_u}$ and $\rho_l=\rho \in (0,1)$,
 we have
 \begin{align}
 \frac{1}{T}\sum_{t=1}^T\|\mathcal{G}(x_t,\nabla F(x_t),\gamma)\|^2  \leq \frac{8\big(\Phi(x_1) -\Phi^* +g(x_1,y_1)-G(x_1)\big)}{T\rho\gamma},
\end{align}
where $\Phi(x)=F(x)+h(x)$.
\end{theorem}

\begin{remark}
Without loss of generality, let $\rho=\rho_l=O(1)$ and $\gamma=\min\big(\frac{3\rho}{4L_F},\frac{\rho\lambda\mu^2}{8\rho_u\hat{L}^2},\frac{\lambda\mu\rho_l}{16L_G\rho_u}, \frac{\mu\rho_l}{16L^2_g\rho_u} \big)=O(1)$.
Since $ F(x_1) -F^* +g(x_1,y_1)-G(x_1)=O(1)$, we have $\frac{1}{T}\sum_{t=1}^T\|\mathcal{G}(x_t,\nabla F(x_t),\gamma)\|^2 \leq O(\frac{1}{T}) \leq \epsilon$. Thus our AdaPAG algorithm obtains the best known convergence rate $O(\frac{1}{T})$, and the best known gradient complexity of $1*T=O(\epsilon^{-1})$ in finding $\epsilon$-stationary solution of Problem~(\ref{eq:1}).
\end{remark}

\subsubsection{Convergence Properties on multimodal $g(x,y)$}
In this subsection, we study the convergence properties of our AdaPAG method when $g(x,\cdot)$ satisfies the local PL condition for all $x$, i.e, it has \textbf{multi local minimizers } $y^\diamond(x,y)\in \arg\min_y f(x,y)$. As~\cite{liu2022bome}, we first define the attraction points.

\begin{definition}
Given any $(x,y)$, if the sequence $\{y_t\}_{t=0}^{\infty}$ generated by gradient descent $y_t=y_{t-1}-\lambda B_t^{-1}\nabla_y g(x,y_{t-1})$ starting from $y_0=y$ converges to $y^\diamond(x,y)$, we say that $y^\diamond(x,y)$ is
the attraction point of $(x,y)$ with step size $\lambda>0$ and adaptive matrix $B_t\succ0$.
\end{definition}

An attraction basin be formed by the same attraction point in set of $(x,y)$. In the following analysis,
we assume the PL condition within the individual attraction basins. Let $F^\diamond(x)=f(x,y^\diamond(x,y))$.

\begin{assumption} \label{ass:1g}
(\textbf{Local PL Condition in Attraction Basins})
Assume that for any $(x,y)$, $y^\diamond(x,y)$ exists.
 $g(x,\cdot)$ satisfies the local PL condition in attraction basins,
i.e., for any $(x,y)$, there exists a constant $\mu>0$ such that
\begin{align}
 \|\nabla_y g(x,y)\|^2 \geq 2\mu \big(g(x,y)- G^\diamond(x)\big),
\end{align}
where $G^\diamond(x)=g(x,y^\diamond(x,y))$.
\end{assumption}

\begin{assumption} \label{ass:2g}
 The function $g\big(x,y^\diamond(x,y)\big)$ satisfies
\begin{align}
 \varrho\big(\nabla^2_{yy}g\big(x,y^\diamond(x,y)\big)\big) \in [\mu, L_g],
\end{align}
where $y^\diamond(x,y)$ is the attraction point of $(x,y)$, and $\varrho(\cdot)$ denotes the eigenvalue (or singular-value) function and $L_g\geq\mu>0$.
\end{assumption}

\begin{assumption} \label{ass:3g}
The functions $f(x,y)$ and $g(x,y)$ satisfy
\begin{itemize}
 \item[(i)] For the attraction points $y^\diamond(x,y)$, we have
 \begin{align}
 {\color{blue}{\|\nabla_yf(x,y^\diamond(x,y))\|\leq C_{fy}, \ \|\nabla^2_{xy}g(x,y^\diamond(x,y))\|\leq C_{gxy}}};
 \end{align}
 \item[(ii)] The partial derivatives $\nabla_x f(x,y)$ and $\nabla_y f(x,y)$ are $L_f$-Lipschitz continuous;
 \item[(iii)] The partial derivatives $\nabla_xg(x,y)$ and $\nabla_yg(x,y)$ are $L_g$-Lipschitz continuous.
\end{itemize}
\end{assumption}

\begin{assumption} \label{ass:5g}
 The function $\Phi^\diamond(x)=F^\diamond(x)+h(x)$ is bounded below in $\mathbb{R}^d$, \emph{i.e.,} $\Phi^\diamond = \inf_{x\in \mathbb{R}^d} \Phi^\diamond(x) > -\infty$.
\end{assumption}

\begin{theorem}  \label{th:1g}
 Assume the sequence $\{x_t,y_t\}_{t=1}^T$ be generated from our Algorithm \ref{alg:1}. Under the above Assumptions~\ref{ass:4}, \ref{ass:6}, \ref{ass:1g}-\ref{ass:5g}, let $0< \gamma \leq \min\big(\frac{3\rho}{4L_F},\frac{\rho\lambda\mu^2}{8\rho_u\hat{L}^2},\frac{\lambda\mu\rho_l}{16L_G\rho_u}, \frac{\mu\rho_l}{16L^2_g\rho_u} \big)$, $0<\lambda\leq \frac{1}{2L_g\rho_u}$ and $\rho_l=\rho \in (0,1)$,
 we have
 \begin{align}
 \frac{1}{T}\sum_{t=1}^T\|\mathcal{G}(x_t,\nabla F^\diamond(x_t),\gamma)\|^2  \leq \frac{8\big(\Phi^\diamond(x_1) -\Phi^\diamond +g(x_1,y_1)-G(x_1)\big)}{T\rho\gamma},
\end{align}
where $\Phi^\diamond(x)=F^\diamond(x)+h(x)$.
\end{theorem}

\begin{remark}
The proof of Theorem~\ref{th:1g} \textbf{can follow} the proof of Theorem~\ref{th:1}.
Our AdaPAG algorithm can also obtain the best known sample (gradient) complexity of $1*T=O(\epsilon^{-1})$ in finding $\epsilon$-stationary solution of Problem~(\ref{eq:1}) \textbf{under local PL condition}.
\end{remark}

\subsection{Convergence Analysis of AdaVSPAG Algorithm}
In this subsection,
we provide the convergence properties of our AdaVSPAG algorithm.
We first provide some mild assumptions for Problem~(\ref{eq:2}).

\begin{assumption} \label{ass:7}
The functions $f(x,y)=\mathbb{E}_{\xi}[f(x,y;\xi)]$, $g(x,y)=\mathbb{E}_{\zeta}[g(x,y;\zeta)]$, $f(x,y;\xi)$ and $g(x,y;\zeta)$ satisfy
\begin{itemize}
 \item[(i)] For the minimizers $y^*(x) \in \arg\min_y g(x,y)$, we have
 \begin{align}
 {\color{blue}{\|\nabla_yf(x,y^*(x))\|\leq C_{fy}, \ \|\nabla^2_{xy}g(x,y^*(x))\|\leq C_{gxy}}};
 \end{align}
 \item[(ii)] The partial derivatives $\nabla_x f(x,y;\xi)$ and $\nabla_y f(x,y;\xi)$ are $L_f$-Lipschitz continuous;
 \item[(iii)] The partial derivatives $\nabla_xg(x,y;\zeta)$ and $\nabla_yg(x,y;\zeta)$ are $L_g$-Lipschitz continuous;
\end{itemize}
\end{assumption}

\begin{assumption} \label{ass:8}
 The partial derivatives $\nabla^2_{xy}g(x,y;\zeta)$ and $\nabla^2_{yy}g(x,y;\zeta)$ are $L_{gxy}$-Lipschitz and $L_{gyy}$-Lipschitz, e.g.,
 for all $x,x_1,x_2 \in \mathbb{R}^d$ and $y,y_1,y_2 \in \mathbb{R}^p$
 \begin{align}
 \|\nabla^2_{xy} g(x_1,y;\zeta)-\nabla^2_{xy} g(x_2,y;\zeta)\| \leq L_{gxy}\|x_1-x_2\|,  \  \|\nabla^2_{xy} g(x,y_1;\zeta)-\nabla^2_{xy} g(x,y_2;\zeta)\| \leq L_{gxy}\|y_1-y_2\|. \nonumber
 \end{align}
 \end{assumption}

\begin{assumption} \label{ass:9}
Stochastic partial derivatives $\nabla_x f(x,y;\xi)$, $\nabla_y f(x,y;\xi)$, $\nabla_x g(x,y;\zeta)$, $\nabla_y g(x,y;\zeta)$, $\nabla^2_{xy} g(x,y;\zeta)$ and $\nabla^2_{yy} g(x,y;\zeta)$ are unbiased
with bounded variance, e.g.,
\begin{align}
\mathbb{E}[\nabla_x f(x,y;\xi)] = \nabla_x f(x,y), \ \mathbb{E}\|\nabla_x f(x,y;\xi) - \nabla_x f(x,y) \|^2 \leq \sigma^2. \nonumber
\end{align}
\end{assumption}

The above Assumptions~\ref{ass:7}-\ref{ass:9} are commonly used in the stochastic bilevel optimization problems \citep{ghadimi2018approximation,yang2021provably,khanduri2021near}.
From Assumptions~\ref{ass:7} and~\ref{ass:8}, we have $\|\nabla_x f(x_1,y)-\nabla_x f(x_2,y)\|=\|\mathbb{E}[\nabla_x f(x_1,y;\xi)-\nabla_x f(x_2,y;\xi)]\|
\leq \mathbb{E}\|\nabla_x f(x_1,y;\xi)-\nabla_x f(x_2,y;\xi)\| \leq L_f\|x_1-x_2\|$, i.e., we can also obtain $\nabla_x f(x,y)$ is $L_f$-Lipschitz continuous,
which is similar for $\nabla_y f(x,y)$, $\nabla_y g(x,y)$, $\nabla^2_{xy}g(x,y)$ and $\nabla^2_{yy}g(x,y)$.
Thus, Assumptions~\ref{ass:7} and~\ref{ass:8} imply that the above Assumptions~\ref{ass:3} and~\ref{ass:4} hold.
Next, we provide a useful lemma.

\begin{lemma} \label{lem:5}
Assume the stochastic gradient $w_t = u_t - G_t(H_t)^{-1}h_t$ be generated from Algorithm \ref{alg:2}, we have
\begin{align}
 \mathbb{E}\|w_t -\nabla F(x_t)\|^2\leq \frac{4\hat{L}^2}{b} \sum_{i=(n_t-1)q}^{t-1}\big( \mathbb{E}\|x_{i+1}-x_i\|^2 + \mathbb{E}\|y_{i+1}-y_i\|^2\big)+  \frac{4\hat{L}^2}{\mu}\big(g(x_t,y_t)-G(x_t)\big) + \frac{\breve{L}^2\sigma^2}{n}, \nonumber
\end{align}
where $\breve{L}^2=8+\frac{8C^2_{fy}}{\mu^2}+ 8\kappa^2+ \frac{8\kappa^2 C^2_{fy}}{\mu^2}$.
\end{lemma}

\begin{theorem} \label{th:3}
Suppose the sequence $\{x_t,y_t\}_{t=1}^T$ be generated from Algorithm \ref{alg:2}. Under the above Assumptions
~\ref{ass:1}, \ref{ass:2}, \ref{ass:5}, \ref{ass:6}, \ref{ass:7}-\ref{ass:9},
let $b=q$, $\gamma\leq \min\Big\{\frac{3\rho}{4L_F},\frac{4\rho\mu^2\lambda}{19\hat{L}^2\rho_u},\frac{\rho }{2\sqrt{38}\hat{L}^2},\frac{\rho}{38\rho_u\lambda \hat{L}^2},\frac{\lambda \rho \rho_u}{4},\frac{\lambda\mu\rho_l}{16L_G\rho_u}, \frac{\mu\rho_l}{16L^2_g\rho_u} \Big\}$, $\rho_l=\rho \in (0,1)$ and $0<\lambda \leq \min\big\{\frac{1}{2L_g\rho_u},\frac{\sqrt{\rho_l}}{4\sqrt{\rho_u}L_g}\big\}$,
we have
\begin{align}
\frac{1}{T}\sum_{t=0}^{T-1}\mathbb{E}\|\mathcal{G}(x_t,\nabla F(x_t),\gamma)\|^2  \leq \frac{32(\Phi(x_0) - \Phi^* + g(x_1,y_1)-G(x_1)) }{3 T\gamma\rho} + \frac{38\breve{L}^2\sigma^2}{3\rho^2 n} +  \frac{32\lambda\sigma^2}{3\rho^2 \gamma n},
\end{align}
where $\breve{L}^2=8+\frac{8C^2_{fy}}{\mu^2}+ 8\kappa^2+ \frac{8\kappa^2 C^2_{fy}}{\mu^2}$ and $\Phi(x)=F(x)+h(x)$.
\end{theorem}

\begin{remark}
Without loss of generality, let $\rho_u=O(1)$, $\rho=\rho_l=O(1)$, $\gamma=O(1)$. Then we have  $\frac{1}{T}\sum_{t=1}^T\|\mathcal{G}(x_t,\nabla F(x_t),\gamma)\|^2 \leq O(\frac{1}{T}) + O(\frac{1}{n})\leq \epsilon$. Meanwhile, let $b=q=\epsilon^{-1/2}$ and $n=T=\epsilon^{-1}$, our AdaVSPAG algorithm can obtain the best known convergence rate $O(\frac{1}{T})$, and the best known sample (gradient) complexity of $b*T + n*T/q=O(\epsilon^{-3/2})$ for finding $\epsilon$-stationary solution of Problem~(\ref{eq:2}).
\end{remark}

Finally, we provide \textbf{an interesting analysis} that our algorithms and convergence results still be competent to the nonconvex strongly-convex bilevel optimization.

\begin{assumption} \label{ass:a1}
 The function $g(x,\cdot)$ is $\mu$-strongly convex on any $x\in \mathbb{R}^d$,
 such that
\begin{align} \label{eq:4.2}
 g(x,y) \geq g(x,y') + \langle\nabla_y g(x,y'), y-y'\rangle + \frac{\mu}{2}\|y-y'\|^2, \ \forall y, y' \in \mathbb{R}^p.
\end{align}
\end{assumption}

By minimizing the inequality~(\ref{eq:4.2}) with respect to $y$, we have
\begin{align}
 \min_y g(x,y) \geq \min_y \Big\{ g(x,y') + \langle\nabla_y g(x,y'), y-y'\rangle + \frac{\mu}{2}\|y-y'\|^2 \Big\}. \label{eq:4.3}
\end{align}
For the right hand side of the inequality~(\ref{eq:4.3}), we have
\begin{align}
 \nabla_y g(x,y') + \mu (y-y') = 0  \Rightarrow y = y' - \frac{1}{\mu}\nabla_y g(x,y').
\end{align}
Then by putting $y = y' - \frac{1}{\mu}\nabla_y g(x,y')$ into the right hand side of the above inequality~(\ref{eq:4.3}), we have
\begin{align}
 \min_y g(x,y) \geq g(x,y') -  \frac{1}{2\mu}\|\nabla_y g(x,y')\|^2. \label{eq:4.4}
\end{align}
Then we can get
\begin{align}
 \|\nabla_y g(x,y')\|^2 \geq 2\mu \big(g(x,y') - \min_y g(x,y)\big) \quad \forall y \in \mathbb{R}^p.
\end{align}
From the above analysis, Assumption~\ref{ass:a1} implies that Assumption~\ref{ass:1} holds.
In other words, Assumption~\ref{ass:1} is milder than Assumption~\ref{ass:a1}.  Since our algorithms and convergence results build on Assumption~\ref{ass:1} instead of Assumption~\ref{ass:a1},
 our algorithms and convergence results still be competent to the nonconvex strongly-convex bilevel optimization.

 Compared with the existing methods for solving nonconvex strongly-convex bilevel optimization, our methods not only use the milder Assumptions~\ref{ass:1} and~\ref{ass:2} than the Assumptions~\ref{ass:a1} used in these existing methods~\citep{ghadimi2018approximation,ji2021bilevel,yang2021provably,khanduri2021near}, but also relax the bounded $\nabla_yf(x,y)$ and $\nabla^2_{xy}g(x,y)$ for all $y\in \mathbb{R}^d$ required in these existing methods~\citep{ghadimi2018approximation,ji2021bilevel,yang2021provably,khanduri2021near}, i.e., our methods only need the bounded $\nabla_yf(x,y^*(x))$ and $\nabla^2_{xy}g(x,y^*(x))$.

\section{Conclusions }
In the paper, we studied a class of nonconvex bilevel optimization, where its upper-level objective is nonconvex possibly with nonsmooth regularization,, and its lower-level objective is nonconvex while satisfies PL condition. We proposed a class of efficient adaptive gradient methods to solve these nonconvex-PL bilevel problems based on mirror descent.
Moreover, we studied the convergence properties of our methods, and proved that our methods obtain the best known convergence rate $O(\frac{1}{T})$.

\small

\bibliographystyle{plainnat}

\bibliography{Ada-NPLBiO}

\newpage

\appendix

\section{Appendix}
In this section, we provide the detailed convergence analysis of our algorithms.
We first review some useful lemmas.

\begin{lemma} \label{lem:A1}
(\cite{nesterov2018lectures})
Assume that $f(x)$ is a differentiable convex function and $\mathcal{X}$ is a convex set.
 $x^* \in \mathcal{X}$ is the solution of the
constrained problem $\min_{x\in \mathcal{X}}f(x)$, if
\begin{align}
 \langle \nabla f(x^*), x-x^*\rangle \geq 0, \ \forall x\in \mathcal{X}.
\end{align}
\end{lemma}

\begin{lemma} \label{lem:A2}
(\cite{karimi2016linear})
 The function $f(x): \mathbb{R}^d\rightarrow \mathbb{R}$ is $L$-smooth and satisfies PL condition with constant $\mu$, then it also satisfies error bound (EB) condition with $\mu$, i.e., for all $x \in \mathbb{R}^d$
\begin{align}
 \|\nabla f(x)\| \geq \mu\|x^*-x\|,
\end{align}
where $x^* \in \arg\min_{x} f(x)$. It also satisfies quadratic growth (QG) condition with $\mu$, i.e.,
\begin{align}
 f(x)-\min_x f(x) \geq \frac{\mu}{2}\|x^*-x\|^2.
\end{align}
\end{lemma}

\begin{lemma} \label{lem:A3}
(\cite{huang2023momentum})
Given the gradient estimator $w_t = u_t - G_t(H_t)^{-1}h_t$ generated from Algorithms~\ref{alg:1} or \ref{alg:2}, for all $t\geq 1$, we have
\begin{align}
 \|w_t-\nabla F(x_t)\|^2
 & \leq 8\|u_t - \nabla_xf(x_t,y_t)\|^2 + \frac{8C^2_{fy}}{\mu^2}\|G_t- \hat{\Pi}_{C_{gxy}}\big[\nabla^2_{xy}g(x_t,y_t)\big] \|^2   \nonumber \\
 & \quad + 8\kappa^2 \|h_t -\Pi_{C_{fy}}\big[\nabla_yf(x_t,y_t)\big] \|^2 + \frac{8\kappa^2 C^2_{fy}}{\mu^2}\|H_t- \nabla^2_{yy} g(x_t,y_t)\big]\|^2
 + \frac{4\hat{L}^2}{\mu}\big(g(x_t,y_t)-G(x_t)\big), \nonumber
\end{align}
where $\kappa=\frac{C_{gxy}}{\mu}$.
\end{lemma}

\subsection{Convergence Analysis of AdaPAG}
In this subsection, we detail the convergence analysis of our AdaPAG algorithm.
We give some useful lemmas.

\begin{lemma} \label{lem:B1}
Suppose the sequence $\{x_t,y_t\}_{t=1}^T$ be generated from Algorithms~\ref{alg:1} or \ref{alg:2}.
Under the above Assumptions~\ref{ass:1}-\ref{ass:3} and \ref{ass:6}, given $\gamma\leq \min\Big\{\frac{\lambda\mu\rho_l}{16L_G\rho_u}, \frac{\mu\rho_l}{16L^2_g\rho_u} \Big\}$, $0<\rho_l\leq 1$ and $0<\lambda \leq \frac{1}{2L_g\rho_u}$, we have
\begin{align}
g(x_{t+1},y_{t+1}) - G(x_{t+1})
& \leq (1-\frac{\lambda\mu}{2\rho_u}) \big(g(x_t,y_t) -G(x_t)\big) + \frac{\rho_l}{8\gamma}\|x_{t+1}-x_t\|^2  -\frac{1}{4\lambda\rho_u}\|y_{t+1}-y_t\|^2 \nonumber \\
& \qquad + \frac{\lambda}{\rho_l}\|\nabla_y g(x_t,y_t)-v_t\|^2,
\end{align}
where $G(x_t)=g(x_t,y^*(x_t))$ with $y^*(x_t) \in \arg\min_{y}g(x_t,y)$ for all $t\geq 1$.
\end{lemma}

\begin{proof}
Using the Assumption~\ref{ass:3}, i.e., $L_g$-smoothness of $g(x,\cdot)$, such that
\begin{align} \label{eq:EE1}
    g(x_{t+1},y_{t+1}) \leq g(x_{t+1},y_t) + \langle \nabla_y g(x_{t+1},y_t), y_{t+1}-y_t \rangle + \frac{L_g}{2}\|y_{t+1}-y_t\|^2.
\end{align}

Since $\rho_u I_{p} \succeq B_t \succeq \rho_l I_{p} \succ 0$ for any $t\geq 1$ is positive definite, we set $B_t=U_t(U_t)^T$, where $ \sqrt{\rho_u} I_{p} \succeq U_t \succeq \sqrt{\rho_l} I_{p} \succ 0$. Thus, we have $B^{-1}_t=(U^{-1}_t)^{T}U^{-1}_t$, where $ \frac{1}{\sqrt{\rho_l}}I_{p} \succeq U^{-1}_t \succeq \frac{1}{\sqrt{\rho_u}}I_{p} \succ 0$.

Since $y_{t+1} = \arg\min_{y\in \mathbb{R}^p}\big\{ \langle v_t, y\rangle
+ \frac{1}{2\lambda}(y-y_t)^TB_t(y-y_t) \big\} = y_t - \lambda B_t^{-1}v_t$, we can obtain
\begin{align} \label{eq:EE2}
    & \langle \nabla_y g(x_{t+1},y_t), y_{t+1}-y_t \rangle \nonumber \\
    & = -\lambda\langle \nabla_y g(x_{t+1},y_t), B^{-1}_tv_t \rangle = -\lambda\langle U^{-1}_t\nabla_y g(x_{t+1},y_t), U^{-1}_tv_t \rangle \nonumber \\
    & = -\frac{\lambda}{2}\Big( \|U^{-1}_t\nabla_y g(x_{t+1},y_t)\|^2 + \|U^{-1}_tv_t\|^2 - \|U^{-1}_t\nabla_y g(x_{t+1},y_t)-U^{-1}_t\nabla_y g(x_t,y_t) + U^{-1}_t\nabla_y g(x_t,y_t)-U^{-1}_tv_t\|^2 \Big) \nonumber \\
    & \leq -\frac{\lambda}{2\rho_u} \|\nabla_y g(x_{t+1},y_t)\|^2 -\frac{1}{2\lambda\rho_u} \|y_{t+1}-y_t\|^2 + \frac{\lambda L^2_g}{\rho_l} \|x_{t+1}-x_t\|^2 + \frac{\lambda}{\rho_l}\|\nabla_y g(x_t,y_t)-v_t\|^2 \nonumber \\
    & \leq -\frac{\lambda\mu}{\rho_u}\big(g(x_{t+1},y_t)-G(x_{t+1})\big)-\frac{1}{2\lambda\rho_u} \|y_{t+1}-y_t\|^2 + \frac{\lambda L^2_g}{\rho_l} \|x_{t+1}-x_t\|^2 + \frac{\lambda}{\rho_l}\|\nabla_y g(x_t,y_t)-v_t\|^2,
\end{align}
where the last inequality is due to the quadratic growth condition of $\mu$-PL functions, i.e.,
\begin{align}
    \|\nabla_y g(x_{t+1},y_t)\|^2 \geq 2\mu\big( g(x_{t+1},y_t)-\min_{y'}g(x_{t+1},y')\big) = 2\mu\big( g(x_{t+1},y_t) - G(x_{t+1})\big).
\end{align}
Substituting \eqref{eq:EE2} into \eqref{eq:EE1}, we have
\begin{align} \label{eq:EE3}
    g(x_{t+1},y_{t+1})
    & \leq g(x_{t+1},y_t)-\frac{\lambda\mu}{\rho_u}\big(g(x_{t+1},y_t)- G(x_{t+1})\big)-\frac{1}{2\lambda\rho_u} \|y_{t+1}-y_t\|^2 + \frac{\lambda L^2_g}{\rho_l} \|x_{t+1}-x_t\|^2 \nonumber \\
    & \quad + \frac{\lambda}{\rho_l}\|\nabla_y g(x_t,y_t)-v_t\|^2 + \frac{L_g}{2}\|y_{t+1}-y_t\|^2,
\end{align}
then rearranging the terms, we can obtain
\begin{align} \label{eq:EE4}
    g(x_{t+1},y_{t+1})- G(x_{t+1})
    & \leq (1-\frac{\lambda\mu}{\rho_u})\big(g(x_{t+1},y_t)- G(x_{t+1})\big)-\frac{1}{2\lambda\rho_u} \|y_{t+1}-y_t\|^2 + \frac{\lambda L^2_g}{\rho_l} \|x_{t+1}-x_t\|^2 \nonumber \\
    & \quad + \frac{\lambda}{\rho_l}\|\nabla_y g(x_t,y_t)-v_t\|^2 + \frac{L_g}{2}\|y_{t+1}-y_t\|^2.
\end{align}

Next, using $L_g$-smoothness of function $f(\cdot,y_t)$, such that
\begin{align}
     g(x_{t+1},y_t) \leq g(x_t,y_t) + \langle \nabla_x g(x_t,y_t), x_{t+1}-x_t \rangle +  \frac{L_g}{2}\|x_{t+1}-x_t\|^2 ,
\end{align}
then we have
\begin{align}
    &g(x_{t+1},y_t) - g(x_t,y_t) \nonumber \\
    & \leq \langle \nabla_x g(x_t,y_t), x_{t+1}-x_t \rangle + \frac{L_g}{2}\|x_{t+1}-x_t\|^2 \nonumber \\
    & = \langle \nabla_x g(x_t,y_t) - \nabla G(x_t), x_{t+1}-x_t \rangle + \langle \nabla G(x_t), x_{t+1}-x_t \rangle + \frac{L_g}{2}\|x_{t+1}-x_t\|^2 \nonumber \\
    & \leq \frac{\rho_l}{8\gamma}\|x_{t+1}-x_t\|^2 + \frac{2\gamma}{\rho_l}\|\nabla_x g(x_t,y_t) - \nabla G(x_t)\|^2  + \langle \nabla G(x_t), x_{t+1}-x_t \rangle + \frac{L_g}{2}\|x_{t+1}-x_t\|^2 \nonumber \\
    & \leq \frac{\rho_l}{8\gamma}\|x_{t+1}-x_t\|^2 + \frac{2L^2_g\gamma}{\rho_l} \|y_t - y^*(x_t)\|^2 + G(x_{t+1}) - G(x_t) \nonumber \\
    & \quad + \frac{L_G}{2}\|x_{t+1}-x_t\|^2 + \frac{L_g}{2}\|x_{t+1}-x_t\|^2 \nonumber \\
    & \leq \frac{4L^2_g\gamma}{\mu\rho_l} \big(g(x_t,y_t) - G(x_t)\big) + G(x_{t+1})- G(x_t)+ (\frac{\rho_l }{8\gamma}+L_G)\|x_{t+1}-x_t\|^2,
\end{align}
where the second last inequality is due to
$L_G$-smoothness of function $G(x)$, and the last inequality holds by Lemma~\ref{lem:A2} and $L_g\leq L_G$.
Then we have
\begin{align} \label{eq:EE5}
    g(x_{t+1},y_t) -G(x_{t+1}) & = g(x_{t+1},y_t)- g(x_t,y_t) + g(x_t,y_t)- G(x_t) + G(x_t) -G(x_{t+1})
    \nonumber \\
    & \leq (1+\frac{4L^2_g\gamma}{\mu\rho_l}) \big(g(x_t,y_t) -G(x_t)\big) + (\frac{\rho_l}{8\gamma}+ L_G)\|x_{t+1}-x_t\|^2.
\end{align}

Substituting \eqref{eq:EE5} in \eqref{eq:EE4}, we get
\begin{align}
    & g(x_{t+1},y_{t+1})- G(x_{t+1})\nonumber \\
    & \leq (1-\frac{\lambda\mu}{\rho_u})(1+\frac{4L^2_g\gamma}{\mu\rho_l}) \big(g(x_t,y_t) -G(x_t)\big) + (1-\frac{\lambda\mu}{\rho_u})(\frac{\rho_l}{8\gamma}+L_G)\|x_{t+1}-x_t\|^2 \nonumber \\
    & \quad -\frac{1}{2\lambda\rho_u} \|y_{t+1}-y_t\|^2 + \frac{\lambda L^2_g}{\rho_l} \|x_{t+1}-x_t\|^2 + \frac{\lambda}{\rho_l}\|\nabla_y g(x_t,y_t)-v_t\|^2 + \frac{L_g}{2}\|y_{t+1}-y_t\|^2 \nonumber \\
    & = (1-\frac{\lambda\mu}{\rho_u})(1+\frac{4L^2_g\gamma}{\mu\rho_l}) \big(g(x_t,y_t) -G(x_t)\big) + \big(\frac{\rho_l}{8\gamma}+L_G-\frac{\lambda\mu}{8\gamma\rho_u}-\frac{L_G\lambda\mu}{\rho_u}
    +\frac{L^2_g\lambda}{\rho_l}\big)\|x_{t+1}-x_t\|^2 \nonumber \\
    & \quad -\frac{1}{2}\big(\frac{1}{\lambda\rho_u}-L_g\big) \|y_{t+1}-y_t\|^2 + \frac{\lambda}{\rho_l}\|\nabla_y g(x_t,y_t)-v_t\|^2 \nonumber \\
    & \leq (1-\frac{\lambda\mu}{2\rho_u}) \big(g(x_t,y_t) -G(x_t)\big) + \frac{\rho_l}{8\gamma}\|x_{t+1}-x_t\|^2  -\frac{1}{4\lambda\rho_u}\|y_{t+1}-y_t\|^2 + \frac{\lambda}{\rho_l}\|\nabla_y g(x_t,y_t)-v_t\|^2,
\end{align}
where the last inequality holds by $\gamma\leq \min\Big\{\frac{\lambda\mu\rho_l}{16L_G\rho_u}, \frac{\mu\rho_l}{16L^2_g\rho_u} \Big\}$, $0<\rho_l\leq 1$, $L_G\geq L_g(1+\kappa)^2$ and $\lambda\leq \frac{1}{2L_g\rho_u}$ for all $t\geq 1$, i.e.,
\begin{align}
   & \gamma\leq \frac{\lambda\mu\rho_l}{16L_G\rho_u} \Rightarrow \lambda \geq \frac{16L_G\gamma\rho_u}{\mu\rho_l} \geq   \frac{16L_g}{\mu\rho_l}(1+\kappa)^2\gamma\rho_u \geq \frac{8\kappa^2\gamma\rho_u}{\rho_l} \Rightarrow  \frac{\lambda\mu}{2\rho_u} \geq \frac{4L^2_g\gamma}{\mu\rho_l} \nonumber \\
   & \gamma \leq \min\Big\{\frac{\lambda\mu\rho_l}{16L_G\rho_u}, \frac{\mu\rho_l}{16L^2_g\rho_u} \Big\}\leq \min\Big\{\frac{\lambda\mu}{16L_G\rho_u}, \frac{\mu\rho_l}{16L^2_g\rho_u} \Big\}
   \Rightarrow \frac{\lambda\mu}{8\gamma\rho_u} \geq
    L_G+\frac{L^2_g\lambda}{\rho_l} \nonumber \\
   &\lambda \leq \frac{1}{2L_g\rho_u}  \Rightarrow \frac{1}{2\lambda\rho_u} \geq  L_g, \ \forall t\geq 1.
\end{align}

\end{proof}

\begin{theorem}  \label{th:A1}
(Restatement of Theorem 1)
 Assume the sequence $\{x_t,y_t\}_{t=1}^T$ be generated from our Algorithm \ref{alg:1}. Under the above Assumptions~\ref{ass:1}-\ref{ass:6}, let $0< \gamma \leq \min\big(\frac{3\rho}{4L_F},\frac{\rho\lambda\mu^2}{8\rho_u\hat{L}^2},\frac{\lambda\mu\rho_l}{16L_G\rho_u}, \frac{\mu\rho_l}{16L^2_g\rho_u} \big)$, $0<\lambda\leq \frac{1}{2L_g\rho_u}$ and $\rho_l=\rho \in (0,1)$,
 we have
 \begin{align}
 \frac{1}{T}\sum_{t=1}^T\|\mathcal{G}(x_t,\nabla F(x_t),\gamma)\|^2  \leq \frac{8\big(\Phi(x_1) -\Phi^* +g(x_1,y_1)-G(x_1)\big)}{T\rho\gamma},
\end{align}
where $\Phi(x)=F(x)+h(x)$.
\end{theorem}

\begin{proof}
By the line 7 of Algorithm~\ref{alg:1}, we have
\begin{align} \label{eq:D1}
x_{t+1} & = \arg\min_{x\in \mathbb{R}^d}\big\{ \langle w_t, x\rangle
+ \frac{1}{2\gamma}(x-x_t)^TA_t(x-x_t) + h(x)\big\} \\
& = \arg\min_{x\in \mathbb{R}^d}\big\{ \frac{1}{2\gamma} \|x-(x_t- A_t^{-1}w_t)\|^2_{A_t} + h(x)\big\} = \mathbb{P}_{h(\cdot)}^\gamma(x_t - A_t^{-1}w_t),
\end{align}
where $\|z\|^2_{A_t} = z^TA_tz$ for all $z\in \mathbb{R}^d$.

Then we define a gradient mapping $\mathcal{G}(x_t,w_t,\gamma) = \frac{1}{\gamma}(x_t-x_{t+1})$.
By the optimality condition of the subproblem~(\ref{eq:D1}), we have for any $z\in \mathbb{R}^d$
\begin{align}
 \big\langle w_t + \frac{1}{\gamma}A_t(x_{t+1}-x_t) + \vartheta_{t+1}, z-x_{t+1}\big\rangle \geq 0,
\end{align}
where $\vartheta_{t+1}\in \partial h(x_{t+1})$.

Let $z=x_t$, and by the convexity of $h(x)$, we can obtain
\begin{align}  \label{eq:D2}
 \langle w_t, x_t - x_{t+1}\rangle & \geq \big\langle\frac{1}{\gamma}A_t(x_{t+1}-x_t), x_{t+1}-x_t\big\rangle + \langle \vartheta_{t+1}, x_{t+1}-x_t\rangle \nonumber \\
 & \geq \frac{\rho}{\gamma}\|x_{t+1}-x_t\|^2 + h(x_{t+1})-h(x_t),
\end{align}
where the last inequality holds by Assumption~\ref{ass:6}, i.e., $A_t\succeq \rho I_d$.

According to the above Lemma \ref{lem:2}, the function $F(x)$ has $L_F$-Lipschitz continuous gradient.
Let $\mathcal{G}(x_t,w_t,\gamma) = \frac{1}{\gamma}(x_t-x_{t+1})$, we have
\begin{align} \label{eq:D3}
  F(x_{t+1}) & \leq F(x_t) + \langle \nabla F(x_t), x_{t+1}-x_t\rangle + \frac{L_F}{2}\|x_{t+1}-x_t\|^2 \nonumber \\
  & = F(x_t) - \gamma \langle \nabla F(x_t),\mathcal{G}(x_t,w_t,\gamma)\rangle + \frac{\gamma ^2L_F}{2}\|\mathcal{G}(x_t,w_t,\gamma)\|^2 \nonumber \\
  & = F(x_t) - \gamma \langle w_t, \mathcal{G}(x_t,w_t,\gamma)\rangle + \gamma \langle w_t - \nabla F(x_t), \mathcal{G}(x_t,w_t,\gamma)\rangle+ \frac{\gamma ^2L_F}{2}\|\mathcal{G}(x_t,w_t,\gamma)\|^2 \nonumber \\
  & \leq F(x_t) - \gamma \rho\|\mathcal{G}(x_t,w_t,\gamma)\|^2 - h(x_{t+1}) + h(x_t) + \gamma \langle w_t - \nabla F(x_t), \mathcal{G}(x_t,w_t,\gamma)\rangle + \frac{\gamma ^2L_F}{2}\|\mathcal{G}(x_t,w_t,\gamma)\|^2 \nonumber \\
  & \leq F(x_t) + ( \frac{\gamma ^2L_F}{2} - \frac{3\gamma \rho}{4})\|\mathcal{G}(x_t,w_t,\gamma)\|^2 - h(x_{t+1}) + h(x_t) + \frac{\gamma }{\rho}\|w_t - \nabla F(x_t)\|^2,
\end{align}
where the second last inequality holds by the above inequality~(\ref{eq:D2}), and the last inequality holds by the following inequality
\begin{align}
 \langle w_t - \nabla F(x_t), \mathcal{G}(x_t,w_t,\gamma)\rangle &\leq \|w_t - \nabla F(x_t)\|\|\mathcal{G}(x_t,w_t,\gamma)\| \nonumber \\
 & \leq  \frac{1}{\rho}\|w_t - \nabla F(x_t)\|^2 + \frac{\rho}{4}\|\mathcal{G}(x_t,w_t,\gamma)\|^2.
\end{align}

Let $\Phi(x) = F(x) + h(x)$. Since $0<\gamma \leq \frac{3\rho}{4L_F}$,
based on the above inequality, we have
\begin{align} \label{eq:D4}
  \Phi(x_{t+1})
  & \leq \Phi(x_t) - \frac{3\gamma \rho}{8}\|\mathcal{G}(x_t,w_t,\gamma)\|^2 + \frac{\gamma }{\rho}\|w_t - \nabla F(x_t)\|^2,
\end{align}

Since $w_t=\hat{\nabla}f(x_t,y_t)$ in Algorithm~\ref{alg:1}, by using Lemma~\ref{lem:3}, we have
\begin{align}\label{eq:D5}
 \|w_t - \nabla F(x_t)\|^2
  = \| \hat{\nabla} f(x_t,y_t) - \nabla F(x_t) \|^2 \leq \frac{2\hat{L}^2}{\mu}\big(g(x_t,y_t)-G(x_t)\big).
\end{align}

Plugging the above inequalities~(\ref{eq:D5}) into (\ref{eq:D4}), we have
\begin{align} \label{eq:D6}
  \Phi(x_{t+1})  \leq \Phi(x_t) +  \frac{2\gamma\hat{L}^2}{\rho\mu}\big(g(x_t,y_t)-G(x_t)\big) - \frac{3\gamma \rho}{8}\|\mathcal{G}(x_t,w_t,\gamma)\|^2,
\end{align}

According to Algorithm~\ref{alg:1}, we have $v_t=\nabla_y g(x_t,y_t)$.
Then by using Lemma \ref{lem:B1}, we have
\begin{align} \label{eq:D7}
& g(x_{t+1},y_{t+1}) - G(x_{t+1}) - \big(g(x_t,y_t) -G(x_t)\big) \nonumber \\
& \leq  -\frac{\lambda\mu}{2\rho_u} \big(g(x_t,y_t) -G(x_t)\big) + \frac{\rho_l}{8\gamma}\|x_{t+1}-x_t\|^2  -\frac{1}{4\lambda\rho_u}\|y_{t+1}-y_t\|^2 + \frac{\lambda}{\rho_l}\|\nabla_y g(x_t,y_t)-v_t\|^2 \nonumber \\
& = -\frac{\lambda\mu}{2\rho_u} \big(g(x_t,y_t) -G(x_t)\big) + \frac{\rho_l\gamma}{8}\|\mathcal{G}(x_t,w_t,\gamma)\|^2  -\frac{1}{4\lambda\rho_u}\|y_{t+1}-y_t\|^2 \nonumber \\
& \leq -\frac{\lambda\mu}{2\rho_u} \big(g(x_t,y_t) -G(x_t)\big) + \frac{\rho_l\gamma}{8}\|\mathcal{G}(x_t,w_t,\gamma)\|^2,
\end{align}
where the above equality holds by $v_t=\nabla_y g(x_t,y_t)$ and $\mathcal{G}(x_t,w_t,\gamma) = \frac{1}{\gamma}(x_t-x_{t+1})$.

Next, we define a useful Lyapunov function (i.e. potential function), for any $t\geq 1$
\begin{align}
 \Omega_t = \Phi(x_t) + g(x_t,y_t)-G(x_t).
\end{align}

By using the above inequalities~(\ref{eq:D6}) and~(\ref{eq:D7}), we have
\begin{align} \label{eq:D8}
 \Omega_{t+1} - \Omega_t & = \Phi(x_{t+1}) - \Phi(x_t) + g(x_{t+1},y_{t+1})-G(x_{t+1})-\big(g(x_t,y_t)-G(x_t)\big)
 \nonumber \\
 & \leq \frac{2\gamma\hat{L}^2}{\rho\mu}\big(g(x_t,y_t)-G(x_t)\big) - \frac{3\gamma \rho}{8}\|\mathcal{G}(x_t,w_t,\gamma)\|^2 -\frac{\lambda\mu}{2\rho_u} \big(g(x_t,y_t) -G(x_t)\big) + \frac{\rho_l\gamma}{8}\|\mathcal{G}(x_t,w_t,\gamma)\|^2   \nonumber \\
& \leq -\frac{\lambda\mu}{4\rho_u} \big(g(x_t,y_t) -G(x_t)\big) - \frac{\rho\gamma}{4}\|\mathcal{G}(x_t,w_t,\gamma)\|^2 ,
\end{align}
where the last inequality is due to $0<\gamma\leq \frac{\rho\lambda\mu^2}{8\rho_u\hat{L}^2}$ and $\rho=\rho_l$.

Since $\mathcal{G}(x_t,w_t,\gamma) = \frac{1}{\gamma}\big(x_t-\mathbb{P}_{h(\cdot)}^{\gamma}(x_t-\gamma A_t^{-1}w_t)\big)$ and $\mathcal{G}(x_t,\nabla F(x_t),\gamma) = \frac{1}{\gamma}\big(x_t-\mathbb{P}_{h(\cdot)}^{\gamma}(x_t-\gamma A_t^{-1}\nabla F(x_t))\big)$,
we have
\begin{align}
 \|\mathcal{G}(x_t,\nabla F(x_t),\gamma)\|^2 & \leq 2\|\mathcal{G}(x_t,w_t,\gamma)\|^2 + 2\|\mathcal{G}(x_t,w_t,\gamma)-\mathcal{G}(x_t,\nabla F(x_t),\gamma)\|^2 \nonumber \\
 & = 2\|\mathcal{G}(x_t,w_t,\gamma)\|^2 + \frac{2}{\gamma^2}\|\mathbb{P}_{h(\cdot)}^{\gamma}(x_t-\gamma A_t^{-1}\nabla F(x_t))-\mathbb{P}_{h(\cdot)}^{\gamma}(x_t-\gamma A_t^{-1}w_t)\|^2 \nonumber \\
 & \leq 2\|\mathcal{G}(x_t,w_t,\gamma)\|^2 + 2\|A_t^{-1}(w_t - \nabla F(x_t))\|^2 \nonumber \\
 & \leq 2\|\mathcal{G}(x_t,w_t,\gamma)\|^2 + 2\|A_t^{-1}\|^2\|w_t - \nabla F(x_t))\|^2 \nonumber \\
 & \leq 2\|\mathcal{G}(x_t,w_t,\gamma)\|^2 + \frac{4\hat{L}^2}{\rho^2\mu}\big(g(x_t,y_t)- G(x_t)\big).
\end{align}
Then we have
\begin{align} \label{eq:D9}
-\|\mathcal{G}(x_t,w_t,\gamma)\|^2 \leq
-\frac{1}{2}\|\mathcal{G}(x_t,\nabla F(x_t),\gamma),\gamma)\|^2 + \frac{2\hat{L}^2}{\rho^2\mu} \big(g(x_t,y_t) -G(x_t)\big).
\end{align}

Plugging the above inequalities~(\ref{eq:D9}) into (\ref{eq:D8}),
we can further get
\begin{align} \label{eq:D10}
\frac{\rho\gamma}{8}\|\mathcal{G}(x_t,\nabla F(x_t),\gamma)\|^2  & \leq \Omega_t  -\Omega_{t+1} + \big(\frac{\hat{L}^2\gamma}{2\rho\mu}-\frac{\lambda\mu}{4\rho_u}\big)\big(g(x_t,y_t)- G(x_t)\big)  \leq \Omega_t  -\Omega_{t+1},
\end{align}
where the last inequality holds by $0<\gamma\leq \frac{\rho\lambda\mu^2}{2\rho_u\hat{L}^2}$.

Based on the above inequality~(\ref{eq:D10}), we have
\begin{align}
\frac{1}{T}\sum_{t=1}^T\|\mathcal{G}(x_t,\nabla F(x_t),\gamma)\|^2 & \leq \frac{1}{T}\sum_{t=1}^T\frac{8(\Omega_t  -\Omega_{t+1})}{\rho\gamma} \nonumber \\
& \leq \frac{8(\Omega_1  - \Phi^*)}{T\rho\gamma} = \frac{8(\Phi(x_1) + g(x_1,y_1)-G(x_1)  - \Phi^*)}{T\rho\gamma},
\end{align}
where the above last inequality holds by Assumption~\ref{ass:5}.

\end{proof}

\subsection{ Convergence Analysis of AdaVSPAG Algorithm }
\label{appendix:C2}
In this subsection, we provide the convergence analysis of our AdaVSPAG algorithm.

\begin{lemma} \label{lem:B2}
Suppose the stochastic gradients $H_t$, $G_t$, $v_t$, $u_t$ and $h_t$ be generated from Algorithm \ref{alg:2}, we have
\begin{align} \label{eq:F1}
 \mathbb{E} \|\mathcal{S}_{[\mu,L_g]}\big[\nabla_{yy}^2 g(x_t,y_t)\big] - H_t\|^2 \leq \frac{2L^2_{gyy}}{b} \sum_{i=(n_t-1)q}^{t-1}\big( \mathbb{E}\|x_{i+1}-x_i\|^2 + \mathbb{E}\|y_{i+1}-y_i\|^2\big) + \frac{\sigma^2}{n},
\end{align}
\begin{align} \label{eq:F2}
 \mathbb{E} \|\hat{\Pi}_{C_{gxy}}\big[\nabla_{xy}^2 g(x_t,y_t)\big] - G_t\|^2 \leq \frac{2L^2_{gxy}}{b} \sum_{i=(n_t-1)q}^{t-1}\big( \mathbb{E}\|x_{i+1}-x_i\|^2 + \mathbb{E}\|y_{i+1}-y_i\|^2\big) + \frac{\sigma^2}{n},
\end{align}
\begin{align} \label{eq:F3}
\mathbb{E} \|\nabla_y g(x_t,y_t) - v_t \|^2 \leq \frac{2L_g^2}{b} \sum_{i=(n_t-1)q}^{t-1}\big( \mathbb{E}\|x_{i+1}-x_i\|^2 + \mathbb{E}\|y_{i+1}-y_i\|^2\big) + \frac{\sigma^2}{n}.
\end{align}
\begin{align} \label{eq:F4}
\mathbb{E} \|\nabla_x f(x_t,y_t) - u_t \|^2 \leq \frac{2L_f^2}{b} \sum_{i=(n_t-1)q}^{t-1}\big( \mathbb{E}\|x_{i+1}-x_i\|^2 + \mathbb{E}\|y_{i+1}-y_i\|^2\big) + \frac{\sigma^2}{n}.
\end{align}
\begin{align} \label{eq:F5}
\mathbb{E} \|\Pi_{C_{fy}}\big[\nabla_y f(x_t,y_t)\big] - h_t \|^2 \leq \frac{2L_f^2}{b} \sum_{i=(n_t-1)q}^{t-1}\big( \mathbb{E}\|x_{i+1}-x_i\|^2 + \mathbb{E}\|y_{i+1}-y_i\|^2\big) + \frac{\sigma^2}{n}.
\end{align}
\end{lemma}

\begin{proof}
We first prove the inequality \eqref{eq:F1}.
According to the definition of $H_t$ in Algorithm \ref{alg:2}, i.e.,
\begin{align}
H_t =\mathcal{S}_{[\mu,L_g]}\big[\frac{1}{b}\sum_{i=1}^b\nabla_{yy}^2 g(x_t,y_t;\zeta^i_t)- \frac{1}{b}\sum_{i=1}^b\nabla_{yy}^2 g (x_{t-1},y_{t-1};\zeta^i_t) + H_{t-1}\big],
\end{align}
we have
\begin{align}
& \mathbb{E} \|\mathcal{S}_{[\mu,L_g]}\big[\nabla_{yy}^2 g(x_t,y_t)\big] - H_t\|^2 \nonumber \\
& \leq \mathbb{E} \|\nabla_{yy}^2 g(x_t,y_t) - \frac{1}{b}\sum_{i=1}^b\nabla_{yy}^2 g(x_t,y_t;\zeta^i_t)+ \frac{1}{b}\sum_{i=1}^b\nabla_{yy}^2 g (x_{t-1},y_{t-1};\zeta^i_t) - H_{t-1}\|^2 \nonumber \\
& = \mathbb{E} \|\nabla_{yy}^2 g(x_t,y_t) - \mathcal{S}_{[\mu,L_g]}\big[\nabla_{yy}^2 g(x_{t-1},y_{t-1})\big] - \frac{1}{b}\sum_{i=1}^b\nabla_{yy}^2 g(x_t,y_t;\zeta^i_t)+ \frac{1}{b}\sum_{i=1}^b\nabla_{yy}^2 g (x_{t-1},y_{t-1};\zeta^i_t) \nonumber \\
& \quad + \mathcal{S}_{[\mu,L_g]}\big[\nabla_{yy}^2 g(x_{t-1},y_{t-1}) \big] - H_{t-1}\|^2 \nonumber \\
& \leq \mathbb{E} \|\nabla_{yy}^2 g(x_t,y_t) - \nabla_{yy}^2 g(x_{t-1},y_{t-1}) - \frac{1}{b}\sum_{i=1}^b\nabla_{yy}^2 g(x_t,y_t;\zeta^i_t)+ \frac{1}{b}\sum_{i=1}^b\nabla_{yy}^2 g (x_{t-1},y_{t-1};\zeta^i_t) \nonumber \\
& \quad + \mathcal{S}_{[\mu,L_g]}\big[\nabla_{yy}^2 g(x_{t-1},y_{t-1}) \big] - H_{t-1}\|^2 \nonumber \\
& = \mathbb{E}\|\mathcal{S}_{[\mu,L_g]}\big[\nabla_{yy}^2 g(x_{t-1},y_{t-1})\big] - H_{t-1}\|^2 + \mathbb{E} \|\nabla_{yy}^2 g(x_t,y_t) - \nabla_{yy}^2 g(x_{t-1},y_{t-1}) - \frac{1}{b}\sum_{i=1}^b\big(\nabla_{yy}^2 g(x_t,y_t;\zeta^i_t) \nonumber \\
& \quad -\nabla_{yy}^2 g (x_{t-1},y_{t-1};\zeta^i_t)\big)\|^2\nonumber \\
& = \mathbb{E}\|\mathcal{S}_{[\mu,L_g]}\big[\nabla_{yy}^2 g(x_{t-1},y_{t-1})\big] - H_{t-1}\|^2 + \frac{1}{b}\mathbb{E} \|\nabla_{yy}^2 g(x_t,y_t) - \nabla_{yy}^2 g(x_{t-1},y_{t-1}) \nonumber \\
& \quad - \big(\nabla_{yy}^2 g(x_t,y_t;\zeta^i_t)-\nabla_{yy}^2 g (x_{t-1},y_{t-1};\zeta^i_t)\big)\|^2 \nonumber \\
& \leq \mathbb{E}\|\mathcal{S}_{[\mu,L_g]}\big[\nabla_{yy}^2 g(x_{t-1},y_{t-1})\big] - H_{t-1}\|^2 + \frac{1}{b}\mathbb{E} \|\nabla_{yy}^2 g(x_t,y_t;\zeta^i_t)-\nabla_{yy}^2 g (x_{t-1},y_{t-1};\zeta^i_t)\|^2 \nonumber \\
& \leq \mathbb{E}\|\mathcal{S}_{[\mu,L_g]}\big[\nabla_{yy}^2 g(x_{t-1},y_{t-1})\big] - H_{t-1}\|^2 + \frac{2L^2_{gyy}}{b} \big( \|x_t-x_{t-1}\|^2 + \|y_t-y_{t-1}\|^2\big), \label{eq:F6}
\end{align}
where the third equality follows by $\mathbb{E}\big[\nabla_{yy}^2 g(x_t,y_t) - \nabla_{yy}^2 g(x_{t-1},y_{t-1}) - \frac{1}{b}\sum_{i=1}^b\big(\nabla_{yy}^2 g(x_t,y_t;\zeta^i_t)-\nabla_{yy}^2 g (x_{t-1},y_{t-1};\zeta^i_t)\big)\big]=0$;
 the second last inequality holds by the inequality $\mathbb{E}\|\zeta - \mathbb{E}[\zeta]\|^2 \leq \mathbb{E}\|\zeta\|^2 $;  the last inequality is due to Assumption~\ref{ass:8}.

Throughout the paper, let $n_t = [t/q]$ such that $(n_t-1)q \leq t \leq n_tq-1$.
Telescoping \eqref{eq:F6} over $t$ from $(n_t-1)q+1$ to $t$, we have
\begin{align}
 \mathbb{E} \|\mathcal{S}_{[\mu,L_g]}\big[\nabla_{yy}^2 g(x_t,y_t)\big] - H_t\|^2 & \leq  \frac{2L^2_{gyy}}{b} \sum_{i=(n_t-1)q}^{t-1}\big( \mathbb{E}\|x_{i+1}-x_i\|^2 + \mathbb{E}\|y_{i+1}-y_i\|^2\big) \nonumber \\
 & \quad + \mathbb{E}\|\mathcal{S}_{[\mu,L_g]}\big[\nabla_{yy}^2 g(x_{(n_t-1)q},y_{(n_t-1)q})\big]- H_{(n_t-1)q}\|^2 \nonumber \\
 & \leq \frac{2L^2_{gyy}}{b} \sum_{i=(n_t-1)q}^{t-1}\big( \mathbb{E}\|x_{i+1}-x_i\|^2 + \mathbb{E}\|y_{i+1}-y_i\|^2\big) + \frac{\sigma^2}{n},
\end{align}
where the last inequality is due to Assumption \ref{ass:9} and $H_{(n_t-1)q}= \mathcal{S}_{[\mu,L_g]}\big[\frac{1}{n}\sum_{i=1}^n \nabla_{yy}^2 g(x_{(n_t-1)q},y_{(n_t-1)q};\zeta_{(n_t-1)q}^i)\big]$.
Similarly, we can obtain
the other above terms.

\end{proof}

\begin{lemma}\label{lem:B3}
Assume the stochastic gradient $w_t = u_t - G_t(H_t)^{-1}h_t$ be generated from Algorithm \ref{alg:2}, we have
\begin{align}
 \mathbb{E}\|w_t -\nabla F(x_t)\|^2\leq \frac{4\hat{L}^2}{b} \sum_{i=(n_t-1)q}^{t-1}\big( \mathbb{E}\|x_{i+1}-x_i\|^2 + \mathbb{E}\|y_{i+1}-y_i\|^2\big)+  \frac{4\hat{L}^2}{\mu}\big(g(x_t,y_t)-G(x_t)\big) + \frac{\breve{L}^2\sigma^2}{n},
\end{align}
where $\breve{L}^2=8+\frac{8C^2_{fy}}{\mu^2}+ 8\kappa^2+ \frac{8\kappa^2 C^2_{fy}}{\mu^2}$.
\end{lemma}

\begin{proof}
According to Lemma~\ref{lem:A3},
we have
\begin{align}
 \|w_t-\nabla F(x_t)\|^2
 & \leq 8\|u_t - \nabla_xf(x_t,y_t)\|^2 + \frac{8C^2_{fy}}{\mu^2}\|G_t-\nabla^2_{xy}g(x_t,y_t)\|^2 + 8\kappa^2 \|h_t -\nabla_yf(x_t,y_t) \|^2  \nonumber \\
 & \quad + \frac{8\kappa^2 C^2_{fy}}{\mu^2}\|H_t-
 \mathcal{S}_{[\mu,L_g]}\big[\nabla^2_{yy} g(x_t,y_t)\big]\|^2 + \frac{4\hat{L}^2}{\mu}\big(g(x_t,y_t)-G(x_t)\big),
\end{align}
where $\hat{L}^2 = 4\big(L^2_f+ \frac{L^2_{gxy}C^2_{fy}}{\mu^2} + \frac{L^2_{gyy} C^2_{gxy}C^2_{fy}}{\mu^4} +
 \frac{L^2_fC^2_{gxy}}{\mu^2}\big)$ and $\kappa=\frac{C_{gxy}}{\mu}$.

Meanwhile, by the above Lemma~\ref{lem:B2}, we can get
\begin{align}
 \mathbb{E}\|w_t-\nabla F(x_t)\|^2 & \leq \frac{4\hat{L}^2}{b} \sum_{i=(n_t-1)q}^{t-1}\big( \mathbb{E}\|x_{i+1}-x_i\|^2 + \mathbb{E}\|y_{i+1}-y_i\|^2\big) + \big(8+\frac{8C^2_{fy}}{\mu^2}+ 8\kappa^2+ \frac{8\kappa^2 C^2_{fy}}{\mu^2}\big)\frac{\sigma^2}{n} \nonumber \\
 & \quad +  \frac{4\hat{L}^2}{\mu}\big(g(x_t,y_t)-G(x_t)\big).
\end{align}

\end{proof}

\begin{theorem} \label{th:A2}
(Restatement of Theorem 3)
Suppose the sequence $\{x_t,y_t\}_{t=1}^T$ be generated from Algorithm \ref{alg:2}.  Under the above Assumptions
~\ref{ass:1}, \ref{ass:2}, \ref{ass:5}, \ref{ass:6}, \ref{ass:7}-\ref{ass:9},
let $b=q$, $\gamma\leq \min\Big\{\frac{3\rho}{4L_F},\frac{4\rho\mu^2\lambda}{19\hat{L}^2\rho_u},\frac{\rho }{2\sqrt{38}\hat{L}^2},\frac{\rho}{38\rho_u\lambda \hat{L}^2},\frac{\lambda \rho \rho_u}{4},\frac{\lambda\mu\rho_l}{16L_G\rho_u}, \frac{\mu\rho_l}{16L^2_g\rho_u} \Big\}$, $\rho_l=\rho \in (0,1)$ and $0<\lambda \leq \min\big\{\frac{1}{2L_g\rho_u},\frac{\sqrt{\rho_l}}{4\sqrt{\rho_u}L_g}\big\}$,
we have
\begin{align}
\frac{1}{T}\sum_{t=0}^{T-1}\mathbb{E}\|\mathcal{G}(x_t,\nabla F(x_t),\gamma)\|^2  \leq \frac{32(\Phi(x_0) - \Phi^* + g(x_1,y_1)-G(x_1)) }{3 T\gamma\rho} + \frac{38\breve{L}^2\sigma^2}{3\rho^2 n} +  \frac{32\lambda\sigma^2}{3\rho^2 \gamma n},
\end{align}
where $\breve{L}^2=8+\frac{8C^2_{fy}}{\mu^2}+ 8\kappa^2+ \frac{8\kappa^2 C^2_{fy}}{\mu^2}$ and $\Phi(x)=F(x)+h(x)$.
\end{theorem}

\begin{proof}
By the line 20 of Algorithm~\ref{alg:2}, we have
\begin{align} \label{eq:H1}
x_{t+1} & = \arg\min_{x\in \mathbb{R}^d}\big\{ \langle w_t, x\rangle
+ \frac{1}{2\gamma}(x-x_t)^TA_t(x-x_t) + h(x)\big\} \\
& = \arg\min_{x\in \mathbb{R}^d}\big\{ \frac{1}{2\gamma} \|x-(x_t- A_t^{-1}w_t)\|^2_{A_t} + h(x)\big\} = \mathbb{P}_{h(\cdot)}^\gamma(x_t - A_t^{-1}w_t),
\end{align}
where $\|z\|^2_{A_t} = z^TA_tz$ for all $z\in \mathbb{R}^d$.

Then we define a gradient mapping $\mathcal{G}(x_t,w_t,\gamma) = \frac{1}{\gamma}(x_t-x_{t+1})$.
By the optimality condition of the subproblem~(\ref{eq:H1}), we have for any $z\in \mathbb{R}^d$
\begin{align}
 \big\langle w_t + \frac{1}{\gamma}A_t(x_{t+1}-x_t) + \vartheta_{t+1}, z-x_{t+1}\big\rangle \geq 0,
\end{align}
where $\vartheta_{t+1}\in \partial h(x_{t+1})$.

Let $z=x_t$, and by the convexity of $h(x)$, we can obtain
\begin{align}  \label{eq:H2}
 \langle w_t, x_t - x_{t+1}\rangle & \geq \big\langle\frac{1}{\gamma}A_t(x_{t+1}-x_t), x_{t+1}-x_t\big\rangle + \langle \vartheta_{t+1}, x_{t+1}-x_t\rangle \nonumber \\
 & \geq \frac{\rho}{\gamma}\|x_{t+1}-x_t\|^2 + h(x_{t+1})-h(x_t),
\end{align}
where the last inequality holds by Assumption~\ref{ass:6}, i.e., $A_t\succeq \rho I_d$.

According to the Lemma \ref{lem:2}, function $F(x)$ has $L_F$-Lipschitz continuous gradient.
Let $\mathcal{G}(x_t,w_t,\gamma) = \frac{1}{\gamma}(x_t-x_{t+1})$, we have
\begin{align} \label{eq:H3}
  F(x_{t+1}) & \leq F(x_t) + \langle \nabla F(x_t), x_{t+1}-x_t\rangle + \frac{L_F}{2}\|x_{t+1}-x_t\|^2 \nonumber \\
  & = F(x_t) - \gamma \langle \nabla F(x_t),\mathcal{G}(x_t,w_t,\gamma)\rangle + \frac{\gamma ^2L_F}{2}\|\mathcal{G}(x_t,w_t,\gamma)\|^2 \nonumber \\
  & = F(x_t) - \gamma \langle w_t, \mathcal{G}(x_t,w_t,\gamma)\rangle + \gamma \langle w_t - \nabla F(x_t), \mathcal{G}(x_t,w_t,\gamma)\rangle+ \frac{\gamma ^2L_F}{2}\|\mathcal{G}(x_t,w_t,\gamma)\|^2 \nonumber \\
  & \leq F(x_t) - \gamma \rho\|\mathcal{G}(x_t,w_t,\gamma)\|^2 - h(x_{t+1}) + h(x_t) + \gamma \langle w_t - \nabla F(x_t), \mathcal{G}(x_t,w_t,\gamma)\rangle + \frac{\gamma ^2L_F}{2}\|\mathcal{G}(x_t,w_t,\gamma)\|^2 \nonumber \\
  & \leq F(x_t) + ( \frac{\gamma ^2L_F}{2} - \frac{3\gamma \rho}{4})\|\mathcal{G}(x_t,w_t,\gamma)\|^2 - h(x_{t+1}) + h(x_t) + \frac{\gamma }{\rho}\|w_t - \nabla F(x_t)\|^2,
\end{align}
where the second last inequality holds by the above inequality~(\ref{eq:H2}), and the last inequality holds by the following inequality
\begin{align}
 \langle w_t - \nabla F(x_t), \mathcal{G}(x_t,w_t,\gamma)\rangle &\leq \|w_t - \nabla F(x_t)\|\|\mathcal{G}(x_t,w_t,\gamma)\| \nonumber \\
 & \leq  \frac{1}{\rho}\|w_t - \nabla F(x_t)\|^2 + \frac{\rho}{4}\|\mathcal{G}(x_t,w_t,\gamma)\|^2.
\end{align}

According to the above Lemma \ref{lem:B3}, we have
\begin{align}\label{eq:H4}
 \mathbb{E}\|w_t - \nabla F(x_t)\|^2
 \leq \frac{4\hat{L}^2}{b} \sum_{i=(n_t-1)q}^{t-1}\big( \mathbb{E}\|x_{i+1}-x_i\|^2 + \mathbb{E}\|y_{i+1}-y_i\|^2\big)+  \frac{4\hat{L}^2}{\mu}\big(g(x_t,y_t)-G(x_t)\big) + \frac{\breve{L}^2\sigma^2}{n}.
\end{align}
Let $\Phi(x) = F(x) + h(x)$. Take expectation into the above inequality~\eqref{eq:H3}, and then plugging~\eqref{eq:H4} into~\eqref{eq:H3}, we have
\begin{align} \label{eq:H5}
 \mathbb{E} [\Phi(x_{t+1})] & \leq \mathbb{E}[\Phi(x_t)] + ( \frac{\gamma ^2L_F}{2} - \frac{3\gamma \rho}{4})\|\mathcal{G}(x_t,w_t,\gamma)\|^2  + \frac{4\gamma\hat{L}^2}{b\rho} \sum_{i=(n_t-1)q}^{t-1}\big( \mathbb{E}\|x_{i+1}-x_i\|^2 + \mathbb{E}\|y_{i+1}-y_i\|^2\big) \nonumber \\
  & \quad +  \frac{4\gamma\hat{L}^2}{\rho\mu}\big(g(x_t,y_t)-G(x_t)\big) + \frac{\gamma\breve{L}^2\sigma^2}{n\rho} \nonumber \\
  & \leq \mathbb{E}[\Phi(x_t)] - \frac{3\gamma \rho}{8} \|\mathcal{G}(x_t,w_t,\gamma)\|^2  + \frac{4\gamma\hat{L}^2}{b\rho} \sum_{i=(n_t-1)q}^{t-1}\big( \mathbb{E}\|x_{i+1}-x_i\|^2 + \mathbb{E}\|y_{i+1}-y_i\|^2\big) \nonumber \\
  & \quad +  \frac{4\gamma\hat{L}^2}{\rho\mu}\big(g(x_t,y_t)-G(x_t)\big) + \frac{\gamma\breve{L}^2\sigma^2}{n\rho}  \nonumber \\
  & =  \mathbb{E}[\Phi(x_t)] - \frac{3\gamma \rho}{16} \|\mathcal{G}(x_t,w_t,\gamma)\|^2 - \frac{3 \rho}{16\gamma} \|x_{t+1}-x_t\|^2 + \frac{4\gamma\hat{L}^2}{b\rho} \sum_{i=(n_t-1)q}^{t-1}\big( \mathbb{E}\|x_{i+1}-x_i\|^2 + \mathbb{E}\|y_{i+1}-y_i\|^2\big) \nonumber \\
  & \quad +  \frac{4\gamma\hat{L}^2}{\rho\mu}\big(g(x_t,y_t)-G(x_t)\big) + \frac{\gamma\breve{L}^2\sigma^2}{n\rho},
\end{align}
where the second inequality is due to $0<\gamma \leq \frac{3\rho}{4L_F}$.

Since $\mathcal{G}(x_t,w_t,\gamma) = \frac{1}{\gamma}\big(x_t-\mathbb{P}_{h(\cdot)}^{\gamma}(x_t-\gamma A_t^{-1}w_t)\big)$ and $\mathcal{G}(x_t,\nabla F(x_t),\gamma) = \frac{1}{\gamma}\big(x_t-\mathbb{P}_{h(\cdot)}^{\gamma}(x_t-\gamma A_t^{-1}\nabla F(x_t))\big)$,
we have
\begin{align}
 \|\mathcal{G}(x_t,\nabla F(x_t),\gamma)\|^2 & \leq 2\|\mathcal{G}(x_t,w_t,\gamma)\|^2 + 2\|\mathcal{G}(x_t,w_t,\gamma)-\mathcal{G}(x_t,\nabla F(x_t),\gamma)\|^2 \nonumber \\
 & = 2\|\mathcal{G}(x_t,w_t,\gamma)\|^2 + \frac{2}{\gamma^2}\|\mathbb{P}_{h(\cdot)}^{\gamma}(x_t-\gamma A_t^{-1}\nabla F(x_t))-\mathbb{P}_{h(\cdot)}^{\gamma}(x_t-\gamma A_t^{-1}w_t)\|^2 \nonumber \\
 & \leq 2\|\mathcal{G}(x_t,w_t,\gamma)\|^2 + 2\|A_t^{-1}(w_t - \nabla F(x_t))\|^2 \nonumber \\
 & \leq 2\|\mathcal{G}(x_t,w_t,\gamma)\|^2 + 2\|A_t^{-1}\|^2\|w_t - \nabla F(x_t))\|^2 \nonumber \\
 & \leq 2\|\mathcal{G}(x_t,w_t,\gamma)\|^2 + \frac{8\hat{L}^2}{b\rho^2} \sum_{i=(n_t-1)q}^{t-1}\big( \mathbb{E}\|x_{i+1}-x_i\|^2 + \mathbb{E}\|y_{i+1}-y_i\|^2\big) \nonumber \\
 & \quad +  \frac{8\hat{L}^2}{\mu\rho^2}\big(g(x_t,y_t)-G(x_t)\big) + \frac{2\breve{L}^2\sigma^2}{n\rho^2}.
\end{align}
Thus we have
\begin{align} \label{eq:H6}
-\|\mathcal{G}(x_t,w_t,\gamma)\|^2 & \leq -\frac{1}{2}\|\mathcal{G}(x_t,\nabla F(x_t),\gamma)\|^2 + \frac{4\hat{L}^2}{b\rho^2} \sum_{i=(n_t-1)q}^{t-1}\big( \mathbb{E}\|x_{i+1}-x_i\|^2 + \mathbb{E}\|y_{i+1}-y_i\|^2\big) \nonumber \\
 & \quad +  \frac{4\hat{L}^2}{\mu\rho^2}\big(g(x_t,y_t)-G(x_t)\big) + \frac{\breve{L}^2\sigma^2}{n\rho^2}.
\end{align}
By plugging \eqref{eq:H6} into \eqref{eq:H5}, we have
\begin{align}
  \mathbb{E} [\Phi(x_{t+1})] & \leq \mathbb{E}[\Phi(x_t)] -\frac{3\gamma\rho}{32}\|\mathcal{G}(x_t,\nabla F(x_t),\gamma)\|^2 - \frac{3 \rho}{16\gamma} \|x_{t+1}-x_t\|^2  + \frac{3\gamma\hat{L}^2}{4b\rho} \sum_{i=(n_t-1)q}^{t-1}\big( \mathbb{E}\|x_{i+1}-x_i\|^2 + \mathbb{E}\|y_{i+1}-y_i\|^2\big) \nonumber \\
 & \quad + \frac{3\gamma\hat{L}^2}{4\mu\rho}\big(g(x_t,y_t)-G(x_t)\big) + \frac{3\gamma\breve{L}^2\sigma^2}{16n\rho}  + \frac{4\gamma\hat{L}^2}{b\rho} \sum_{i=(n_t-1)q}^{t-1}\big( \mathbb{E}\|x_{i+1}-x_i\|^2 + \mathbb{E}\|y_{i+1}-y_i\|^2\big) \nonumber \\
  & \quad +  \frac{4\gamma\hat{L}^2}{\rho\mu}\big(g(x_t,y_t)-G(x_t)\big) + \frac{\gamma\breve{L}^2\sigma^2}{n\rho} \nonumber \\
  & \leq \mathbb{E}[\Phi(x_t)] -\frac{3\gamma\rho}{32}\|\mathcal{G}(x_t,\nabla F(x_t),\gamma)\|^2 - \frac{3 \rho}{16\gamma} \|x_{t+1}-x_t\|^2 + \frac{19\gamma\hat{L}^2}{4b\rho} \sum_{i=(n_t-1)q}^{t-1}\big( \mathbb{E}\|x_{i+1}-x_i\|^2 + \mathbb{E}\|y_{i+1}-y_i\|^2\big) \nonumber \\
 & \quad +  \frac{19\gamma\hat{L}^2}{4\mu\rho}\big(g(x_t,y_t)-G(x_t)\big) + \frac{19\gamma\breve{L}^2\sigma^2}{16n\rho}.
\end{align}

Next, we define a useful Lyapunov function, for any $t\geq 1$
\begin{align}
 \Psi_t = \mathbb{E}\big[ \Phi(x_t) + g(x,y)-G(x) \big].
\end{align}
According to Lemma \ref{lem:B1}, we have
\begin{align}
& g(x_{t+1},y_{t+1}) - G(x_{t+1}) - \big(g(x_t,y_t) -G(x_t)\big) \nonumber \\
& \leq  -\frac{\lambda\mu}{2\rho_u} \big(g(x_t,y_t) -G(x_t)\big) + \frac{\rho_l}{8\gamma}\|x_{t+1}-x_t\|^2  -\frac{1}{4\lambda\rho_u}\|y_{t+1}-y_t\|^2 + \frac{\lambda}{\rho_l}\|\nabla_y g(x_t,y_t)-v_t\|^2.
\end{align}

Then we have
\begin{align} \label{eq:H7}
 \Psi_{t+1} - \Psi_t & = \mathbb{E}\big[ \Phi(x_{t+1}) - \Phi(x_t) + g(x_{t+1},y_{t+1}) - G(x_{t+1}) - \big(g(x_t,y_t) -G(x_t)\big) \big] \nonumber \\
 & \leq  -\frac{3\gamma\rho}{32}\mathbb{E}\|\mathcal{G}(x_t,\nabla F(x_t),\gamma)\|^2 - \frac{3 \rho}{16\gamma} \mathbb{E}\|x_{t+1}-x_t\|^2 + \frac{19\gamma\hat{L}^2}{4b\rho} \sum_{i=(n_t-1)q}^{t-1}\big( \mathbb{E}\|x_{i+1}-x_i\|^2 + \mathbb{E}\|y_{i+1}-y_i\|^2\big) \nonumber \\
 & \quad +  \frac{19\gamma\hat{L}^2}{4\mu\rho}\big(g(x_t,y_t)-G(x_t)\big) + \frac{19\gamma\breve{L}^2\sigma^2}{16n\rho}  \nonumber \\
 & \quad -\frac{\lambda\mu}{2\rho_u} \big(g(x_t,y_t) -G(x_t)\big) + \frac{\rho_l}{8\gamma}\mathbb{E}\|x_{t+1}-x_t\|^2  -\frac{1}{4\lambda\rho_u}\mathbb{E}\|y_{t+1}-y_t\|^2 + \frac{\lambda}{\rho_l}\mathbb{E}\|\nabla_y g(x_t,y_t)-v_t\|^2 \nonumber \\
 & = - \frac{3\gamma \rho}{32} \mathbb{E}\|\mathcal{G}(x_t,\nabla F(x_t),\gamma)\|^2  - \big(\frac{3\rho}{16\gamma}- \frac{\rho_l}{8\gamma}\big)\mathbb{E}\|x_{t+1}-x_t\|^2 + \frac{19\gamma\hat{L}^2}{4b\rho} \sum_{i=(n_t-1)q}^{t-1}\big( \mathbb{E}\|x_{i+1}-x_i\|^2 + \mathbb{E}\|y_{i+1}-y_i\|^2\big) \nonumber \\
 & \quad - \big(\frac{\lambda\mu}{2\rho_u} - \frac{19\gamma \hat{L}^2 }{4\rho\mu} \big) \big(g(x_t,y_t) -G(x_t)\big) -\frac{1}{4\lambda\rho_u}\|y_{t+1}-y_t\|^2 + \frac{\lambda}{\rho_l}\|\nabla_y g(x_t,y_t)-v_t\|^2  + \frac{19\gamma\breve{L}^2\sigma^2}{16n\rho} \nonumber \\
 & \leq - \frac{3\gamma \rho}{32} \mathbb{E}\|\mathcal{G}(x_t,\nabla F(x_t),\gamma)\|^2  - \frac{\rho}{16\gamma}\mathbb{E}\|x_{t+1}-x_t\|^2 + \frac{19\gamma\hat{L}^2}{4b\rho} \sum_{i=(n_t-1)q}^{t-1}\big( \mathbb{E}\|x_{i+1}-x_i\|^2 + \mathbb{E}\|y_{i+1}-y_i\|^2\big) \nonumber \\
 & \quad -\frac{1}{4\lambda\rho_u}\|y_{t+1}-y_t\|^2 + \frac{\lambda}{\rho_l}\mathbb{E}\|\nabla_y g(x_t,y_t)-v_t\|^2 + \frac{19\gamma\breve{L}^2\sigma^2}{16n\rho} ,
\end{align}
where the last inequality holds by $\gamma \leq \frac{4\rho\mu^2\lambda}{19\hat{L}^2\rho_u}$ and $\rho=\rho_l$.

Averaging over $t=0,1,\cdots,T-1$ on both sides of \eqref{eq:H7}, we have
\begin{align} \label{eq:H8}
&\frac{3\gamma\rho}{32}\frac{1}{T}\sum_{t=0}^{T-1}\mathbb{E}\|\mathcal{G}(x_t,\nabla F(x_t),\gamma)\|^2 \nonumber \\
& \leq \frac{\Psi_0 - \Psi_T}{T}  -\frac{\rho}{16\gamma}\frac{1}{T}\sum_{t=0}^{T-1}\mathbb{E}\|x_{t+1}-x_t\|^2 - \frac{1}{4\lambda\rho_u}\frac{1}{T}\sum_{t=0}^{T-1}\|y_{t+1}-y_t\|^2 \nonumber \\
& \quad + \frac{19\gamma}{8\rho} \frac{1}{T}\sum_{t=0}^{T-1}\Big( \frac{2\hat{L}^2}{b} \sum_{i=(n_t-1)q}^{t-1}\big( \mathbb{E}\|x_{i+1}-x_i\|^2 + \mathbb{E}\|y_{i+1}-y_i\|^2\big) + \frac{\breve{L}^2\sigma^2}{2n} \Big)  \nonumber \\
& \quad + \frac{\lambda}{\rho_l} \frac{1}{T}\sum_{t=0}^{T-1}\Big( \frac{2L_g^2}{b} \sum_{i=(n_t-1)q}^{t-1}\big( \mathbb{E}\|x_{i+1}-x_i\|^2 + \mathbb{E}\|y_{i+1}-y_i\|^2\big) + \frac{\sigma^2}{n} \Big)  \nonumber \\
& \leq \frac{\Psi_0 - \Psi_T}{T}  -\frac{\rho}{16\gamma}\frac{1}{T}\sum_{t=0}^{T-1}\mathbb{E}\|x_{t+1}-x_t\|^2 - \frac{1}{4\lambda\rho_u}\frac{1}{T}\sum_{t=0}^{T-1}\|y_{t+1}-y_t\|^2 \nonumber \\
& \quad + \frac{19\gamma}{8\rho}\frac{1}{T}\sum_{t=0}^{T-1}\Big( \frac{2q\hat{L}^2}{b}\big( \mathbb{E}\|x_{t+1}-x_t\|^2 + \mathbb{E}\|y_{t+1}-y_t\|^2\big) + \frac{\breve{L}^2\sigma^2}{2n} \Big)  \nonumber \\
& \quad + \frac{\lambda}{\rho_l}\frac{1}{T}\sum_{t=0}^{T-1}\Big( \frac{2qL_g^2}{b}\big( \mathbb{E}\|x_{t+1}-x_t\|^2 + \mathbb{E}\|y_{t+1}-y_t\|^2\big) + \frac{\sigma^2}{n} \Big)   \nonumber \\
& = \frac{\Psi_0 - \Psi_T}{T}  - \big(\frac{\rho}{16\gamma} -\frac{19\gamma \hat{L}^2q}{4\rho b} -\frac{2\lambda L^2_gq}{\rho_l b} \big)\frac{1}{T}\sum_{t=0}^{T-1}\mathbb{E}\|x_{t+1}-x_t\|^2 \nonumber \\
& \quad - \big(\frac{1}{4\lambda\rho_u} - \frac{19\gamma \hat{L}^2q}{4\rho b} -\frac{2\lambda L^2_gq}{\rho_l b} \big)\frac{1}{T}\sum_{t=0}^{T-1}\mathbb{E}\|y_{t+1}-y_t\|^2 + \frac{19\gamma \breve{L}^2 \sigma^2}{16\rho n} +  \frac{\lambda\sigma^2}{\rho_l n},
\end{align}
where the first inequality is due to by lemma \ref{lem:B2}, and the second inequality holds by $\sum_{t=0}^{T-1}\sum_{i=(n_t-1)q}^{t-1}\big( \mathbb{E}\|x_{i+1}-x_i\|^2 + \mathbb{E}\|y_{i+1}-y_i\|^2\big) \leq q\sum_{t=0}^{T-1}\big( \mathbb{E}\|x_{t+1}-x_t\|^2 + \mathbb{E}\|y_{t+1}-y_t\|^2\big)$.

Let $b=q$, $0<\lambda\leq \frac{\sqrt{\rho_l}}{4\sqrt{\rho_u}L_g}$, $0< \gamma \leq \min\big(\frac{\rho }{2\sqrt{38}\hat{L}^2},\frac{\rho}{38\rho_u\lambda \hat{L}^2},\frac{\lambda \rho \rho_u}{4}\big)$, we have
\begin{align}
& 0< \gamma \leq \min\big(\frac{\rho }{2\sqrt{38}\hat{L}^2},\frac{\rho}{38\rho_u\lambda \hat{L}^2}\big) \Rightarrow \frac{\rho}{16\gamma} -\frac{19\gamma \hat{L}^2q}{4\rho b} -\frac{2\lambda L^2_gq}{\rho_l b} \geq 0 \label{eq:H9} \\
& 0<\lambda\leq \frac{\sqrt{\rho_l}}{4\sqrt{\rho_u}L_g}, \ 0< \gamma \leq \min\big(\frac{\rho }{2\sqrt{38}\hat{L}^2},\frac{\rho}{38\rho_u\lambda \hat{L}^2},\frac{\lambda \rho \rho_u}{4}\big) \Rightarrow
\frac{1}{4\lambda\rho_u} - \frac{19\gamma \hat{L}^2q}{4\rho b} -\frac{2\lambda L^2_gq}{\rho_l b} \geq 0. \label{eq:H10}
\end{align}

Based on the above inequalities \eqref{eq:H8}, \eqref{eq:H9} and \eqref{eq:H10}, we have
\begin{align}
\frac{3\gamma\rho}{32}\frac{1}{T}\sum_{t=0}^{T-1}\mathbb{E}\|\mathcal{G}(x_t,\nabla F(x_t),\gamma)\|^2 & \leq \frac{\Psi_0 - \Psi_T }{T} + \frac{19\gamma \breve{L}^2 \sigma^2}{16\rho n} +  \frac{\lambda\sigma^2}{\rho_l n} \nonumber \\
& \leq \frac{\Phi(x_0) - \Phi^* + g(x_1,y_1)-G(x_1) }{T} + \frac{19\gamma \breve{L}^2 \sigma^2}{16\rho n} +  \frac{\lambda\sigma^2}{\rho_l n},
\end{align}
where the last inequality holds by Assumption~\ref{ass:5}.

Then we can obtain
\begin{align}
\frac{1}{T}\sum_{t=0}^{T-1}\mathbb{E}\|\mathcal{G}(x_t,\nabla F(x_t),\gamma)\|^2  \leq \frac{32(\Phi(x_0) - \Phi^* + g(x_1,y_1)-G(x_1)) }{3 T\gamma\rho} + \frac{38\breve{L}^2\sigma^2}{3\rho^2 n} +  \frac{32\lambda\sigma^2}{3\rho^2 \gamma n}.
\end{align}

\end{proof}

\end{document}